

\magnification 1250
\parindent=0cm   \parskip=0pt
\pageno=1

\def\ind{\par\hskip 1cm\relax}

\pretolerance=500 \tolerance=1000  \brokenpenalty=5000

\catcode`\@=11

\font\eightrm=cmr8         \font\eighti=cmmi8
\font\eightsy=cmsy8        \font\eightbf=cmbx8
\font\eighttt=cmtt8        \font\eightit=cmti8
\font\eightsl=cmsl8        \font\sixrm=cmr6
\font\sixi=cmmi6           \font\sixsy=cmsy6
\font\sixbf=cmbx6


\font\tengoth=eufm10       \font\tenbboard=msbm10
\font\eightgoth=eufm8      \font\eightbboard=msbm8
\font\sevengoth=eufm7      \font\sevenbboard=msbm7
\font\sixgoth=eufm6        \font\fivegoth=eufm5

\font\tencyr=wncyr10       
\font\eightcyr=wncyr8      
\font\sevencyr=wncyr7      
\font\sixcyr=wncyr6


\skewchar\eighti='177 \skewchar\sixi='177
\skewchar\eightsy='60 \skewchar\sixsy='60


\newfam\gothfam           \newfam\bboardfam
\newfam\cyrfam

\def\tenpoint{%
  \textfont0=\tenrm \scriptfont0=\sevenrm \scriptscriptfont0=\fiverm
  \def\rm{\fam\z@\tenrm}%
  \textfont1=\teni  \scriptfont1=\seveni  \scriptscriptfont1=\fivei
  \def\oldstyle{\fam\@ne\teni}\let\old=\oldstyle
  \textfont2=\tensy \scriptfont2=\sevensy \scriptscriptfont2=\fivesy
  \textfont\gothfam=\tengoth \scriptfont\gothfam=\sevengoth
  \scriptscriptfont\gothfam=\fivegoth
  \def\goth{\fam\gothfam\tengoth}%
  \textfont\bboardfam=\tenbboard \scriptfont\bboardfam=\sevenbboard
  \scriptscriptfont\bboardfam=\sevenbboard
  \def\bb{\fam\bboardfam\tenbboard}%
 \textfont\cyrfam=\tencyr \scriptfont\cyrfam=\sevencyr
  \scriptscriptfont\cyrfam=\sixcyr
  \def\cyr{\fam\cyrfam\tencyr}%
  \textfont\itfam=\tenit
  \def\it{\fam\itfam\tenit}%
  \textfont\slfam=\tensl
  \def\sl{\fam\slfam\tensl}%
  \textfont\bffam=\tenbf \scriptfont\bffam=\sevenbf
  \scriptscriptfont\bffam=\fivebf
  \def\bf{\fam\bffam\tenbf}%
  \textfont\ttfam=\tentt
  \def\tt{\fam\ttfam\tentt}%
  \abovedisplayskip=9pt plus 3pt minus 9pt
  \belowdisplayskip=\abovedisplayskip
  \abovedisplayshortskip=0pt plus 3pt
  \belowdisplayshortskip=4pt plus 3pt
  \smallskipamount=3pt plus 1pt minus 1pt
  \medskipamount=6pt plus 2pt minus 2pt
  \bigskipamount=12pt plus 4pt minus 4pt
  \normalbaselineskip=12pt
  \setbox\strutbox=\hbox{\vrule height8.5pt depth3.5pt width0pt}%
  \let\bigf@nt=\tenrm       \let\smallf@nt=\sevenrm
  \normalbaselines\rm}

\def\eightpoint{%
  \textfont0=\eightrm \scriptfont0=\sixrm \scriptscriptfont0=\fiverm
  \def\rm{\fam\z@\eightrm}%
  \textfont1=\eighti  \scriptfont1=\sixi  \scriptscriptfont1=\fivei
  \def\oldstyle{\fam\@ne\eighti}\let\old=\oldstyle
  \textfont2=\eightsy \scriptfont2=\sixsy \scriptscriptfont2=\fivesy
  \textfont\gothfam=\eightgoth \scriptfont\gothfam=\sixgoth
  \scriptscriptfont\gothfam=\fivegoth
  \def\goth{\fam\gothfam\eightgoth}%
  \textfont\cyrfam=\eightcyr \scriptfont\cyrfam=\sixcyr
  \scriptscriptfont\cyrfam=\sixcyr
  \def\cyr{\fam\cyrfam\eightcyr}%
  \textfont\bboardfam=\eightbboard \scriptfont\bboardfam=\sevenbboard
  \scriptscriptfont\bboardfam=\sevenbboard
  \def\bb{\fam\bboardfam}%
  \textfont\itfam=\eightit
  \def\it{\fam\itfam\eightit}%
  \textfont\slfam=\eightsl
  \def\sl{\fam\slfam\eightsl}%
  \textfont\bffam=\eightbf \scriptfont\bffam=\sixbf
  \scriptscriptfont\bffam=\fivebf
  \def\bf{\fam\bffam\eightbf}%
  \textfont\ttfam=\eighttt
  \def\tt{\fam\ttfam\eighttt}%
  \abovedisplayskip=9pt plus 3pt minus 9pt
  \belowdisplayskip=\abovedisplayskip
  \abovedisplayshortskip=0pt plus 3pt
  \belowdisplayshortskip=3pt plus 3pt
  \smallskipamount=2pt plus 1pt minus 1pt
  \medskipamount=4pt plus 2pt minus 1pt
  \bigskipamount=9pt plus 3pt minus 3pt
  \normalbaselineskip=9pt
  \setbox\strutbox=\hbox{\vrule height7pt depth2pt width0pt}%
  \let\bigf@nt=\eightrm     \let\smallf@nt=\sixrm
  \normalbaselines\rm}

\tenpoint
\def\pc#1{\bigf@nt#1\smallf@nt}         \def\pd#1 {{\pc#1} }
\frenchspacing
\def\raggedbottom{\topskip 10pt plus 36pt\r@ggedbottomtrue}

\mathcode`A="7041 \mathcode`B="7042 \mathcode`C="7043 \mathcode`D="7044
\mathcode`E="7045 \mathcode`F="7046 \mathcode`G="7047 \mathcode`H="7048
\mathcode`I="7049 \mathcode`J="704A \mathcode`K="704B \mathcode`L="704C
\mathcode`M="704D \mathcode`N="704E \mathcode`O="704F \mathcode`P="7050
\mathcode`Q="7051 \mathcode`R="7052 \mathcode`S="7053 \mathcode`T="7054
\mathcode`U="7055 \mathcode`V="7056 \mathcode`W="7057 \mathcode`X="7058
\mathcode`Y="7059 \mathcode`Z="705A

\def\spacedmath#1{\def\packedmath##1${\bgroup\mathsurround=0pt ##1\egroup$}%
\mathsurround#1 \everymath={\packedmath}\everydisplay={\mathsurround=0pt }}

\def\nospacedmath{\mathsurround=0pt \everymath={}\everydisplay={} }

\def\diagram#1{\def\normalbaselines{\baselineskip=0truept
\lineskip=10truept\lineskiplimit=1truept}   \matrix{#1}}

\def\hfl#1#2{\normalbaselines{\baselineskip=0truept
\lineskip=10truept\lineskiplimit=1truept}\nospacedmath\smash{\mathop{\hbox to
12truemm{\rightarrowfill}}\limits^{\scriptstyle#1}_{\scriptstyle#2}}}

\def\ghfl#1#2{\normalbaselines{\baselineskip=0pt
\lineskip=10truept\lineskiplimit=1truept}\nospacedmath\smash{\mathop{\hbox to
16truemm{\rightarrowfill}}\limits^{\scriptstyle#1}_{\scriptstyle#2}}}

\def\phfl#1#2{\normalbaselines{\baselineskip=0pt
\lineskip=10truept\lineskiplimit=1truept}\nospacedmath\smash{\mathop{\hbox to
8truemm{\rightarrowfill}}\limits^{\scriptstyle#1}_{\scriptstyle#2}}}

\def\vfl#1#2{\llap{$\scriptstyle#1$}\left\downarrow\vbox to
6truemm{}\right.\rlap{$\scriptstyle#2$}}

\def\up#1{\raise 1ex\hbox{\smallf@nt#1}}
\def\tx{\kern-1.5pt -}
\def\cqfd{\kern 2truemm\unskip\penalty 500\vrule height 4pt depth 0pt width
4pt\medbreak}
\def\carre{\vrule height 4pt depth 0pt width 4pt}
\def\virg{\raise .4ex\hbox{,}}   
\def\decale#1{\smallbreak\hskip 28pt\llap{#1}\kern 5pt}
\def\no{n\up{o}\kern 2pt}

\def\Medbreak{\vskip-\lastskip\medbreak}
\def\pointir{\unskip . --- \ignorespaces}
\long\def\th#1 #2\enonce#3\endth{%
   \Medbreak
   {\pc#1} {#2\unskip}\pointir{\it #3}\medskip}

\def\crog{{\vrule height 2.57mm depth 0.85mm width 0.3mm}\kern -0.3mm [}

\def\crod{]\kern -0.4mm{\vrule height 2.57mm depth 0.85mm
width 0.3 mm}}

\def\moins{\mathrel{\hbox{\vrule height 3pt depth -2pt
width 6pt}}}

\def\rond{\kern 1pt{\scriptstyle\circ}\kern 1pt}
\def\epi{\rightarrow \kern -3mm\rightarrow }
\def\longepi{\longrightarrow \kern -5mm\longrightarrow }
\def\mono{\lhook\joinrel\mathrel{\longrightarrow}}
\def\iso{\mathrel{\mathop{\kern 0pt\longrightarrow }\limits^{\sim}}}
\def\End{\mathop{\rm End}\nolimits}
\def\Hom{\mathop{\rm Hom}\nolimits}

\def\im{\mathop{\rm Im}\nolimits}

\def\dim{\mathop{\rm dim}\nolimits}
\def\Card{\mathop{\rm Card}\nolimits}
\def\Tr{\mathop{\rm Tr}\nolimits}

\def\Res{\mathop{\rm Res}\nolimits}

\def\tx{\kern-1.5pt -}
\def\note#1#2{\footnote{\parindent
0.4cm$^#1$}{\vtop{\eightpoint\baselineskip12pt\hsize15.5truecm\noindent
#2}}\parindent 0cm}

\def\build#1_#2^#3{\mathrel{
\mathop{\kern 0pt#1}\limits_{#2}^{#3}}}

\def\fleche(#1,#2)\dir(#3,#4)\long#5{%
\noalign{\leftput(#1,#2){\vector(#3,#4){#5}}}}

\def\ligne(#1,#2)\dir(#3,#4)\long#5{%
\noalign{\leftput(#1,#2){\lline(#3,#4){#5}}}}

\def\put(#1,#2)#3{\noalign{\setbox1=\hbox{%
    \kern #1\unitlength
    \raise #2\unitlength\hbox{$#3$}}%
    \ht1=0pt \wd1=0pt \dp1=0pt\box1}}
 \message{`lline' & `vector' macros from LaTeX}
 \catcode`@=11
\def\{{\relax\ifmmode\lbrace\else$\lbrace$\fi}
\def\}{\relax\ifmmode\rbrace\else$\rbrace$\fi}
\def\newcount{\alloc@0\count\countdef\insc@unt}
\def\newdimen{\alloc@1\dimen\dimendef\insc@unt}
\def\newwrite{\alloc@7\write\chardef\sixt@@n}

\newwrite\@unused
\newcount\@tempcnta
\newcount\@tempcntb
\newdimen\@tempdima
\newdimen\@tempdimb
\newbox\@tempboxa

\def\@spaces{\space\space\space\space}
\def\@whilenoop#1{}
\def\@whiledim#1\do #2{\ifdim #1\relax#2\@iwhiledim{#1\relax#2}\fi}
\def\@iwhiledim#1{\ifdim #1\let\@nextwhile=\@iwhiledim
        \else\let\@nextwhile=\@whilenoop\fi\@nextwhile{#1}}

\font\tenln    = line10
\font\tenlnw   = linew10
\newdimen\@wholewidth
\newdimen\@halfwidth
\newdimen\unitlength

\unitlength =1pt
\def\thinlines{\let\@linefnt\tenln \let\@circlefnt\tencirc
  \@wholewidth\fontdimen8\tenln \@halfwidth .5\@wholewidth}
\def\thicklines{\let\@linefnt\tenlnw \let\@circlefnt\tencircw
  \@wholewidth\fontdimen8\tenlnw \@halfwidth .5\@wholewidth}

\def\linethickness#1{\@wholewidth #1\relax \@halfwidth .5\@wholewidth}

\newif\if@negarg

\def\lline(#1,#2)#3{\@xarg #1\relax \@yarg #2\relax
\@linelen=#3\unitlength
\ifnum\@xarg =0 \@vline
  \else \ifnum\@yarg =0 \@hline \else \@sline\fi
\fi}

\def\@sline{\ifnum\@xarg< 0 \@negargtrue \@xarg -\@xarg \@yyarg -\@yarg
  \else \@negargfalse \@yyarg \@yarg \fi
\ifnum \@yyarg >0 \@tempcnta\@yyarg \else \@tempcnta -\@yyarg \fi
\ifnum\@tempcnta>6 \@badlinearg\@tempcnta0 \fi
\setbox\@linechar\hbox{\@linefnt\@getlinechar(\@xarg,\@yyarg)}%
\ifnum \@yarg >0 \let\@upordown\raise \@clnht\z@
   \else\let\@upordown\lower \@clnht \ht\@linechar\fi
\@clnwd=\wd\@linechar
\if@negarg \hskip -\wd\@linechar \def\@tempa{\hskip -2\wd\@linechar}\else
     \let\@tempa\relax \fi
\@whiledim \@clnwd <\@linelen \do
  {\@upordown\@clnht\copy\@linechar
   \@tempa
   \advance\@clnht \ht\@linechar
   \advance\@clnwd \wd\@linechar}%
\advance\@clnht -\ht\@linechar
\advance\@clnwd -\wd\@linechar
\@tempdima\@linelen\advance\@tempdima -\@clnwd
\@tempdimb\@tempdima\advance\@tempdimb -\wd\@linechar
\if@negarg \hskip -\@tempdimb \else \hskip \@tempdimb \fi
\multiply\@tempdima \@m
\@tempcnta \@tempdima \@tempdima \wd\@linechar \divide\@tempcnta \@tempdima
\@tempdima \ht\@linechar \multiply\@tempdima \@tempcnta
\divide\@tempdima \@m
\advance\@clnht \@tempdima
\ifdim \@linelen <\wd\@linechar
   \hskip \wd\@linechar
  \else\@upordown\@clnht\copy\@linechar\fi}

\def\@hline{\ifnum \@xarg <0 \hskip -\@linelen \fi
\vrule height \@halfwidth depth \@halfwidth width \@linelen
\ifnum \@xarg <0 \hskip -\@linelen \fi}

\def\@getlinechar(#1,#2){\@tempcnta#1\relax\multiply\@tempcnta 8
\advance\@tempcnta -9 \ifnum #2>0 \advance\@tempcnta #2\relax\else
\advance\@tempcnta -#2\relax\advance\@tempcnta 64 \fi
\char\@tempcnta}

\def\vector(#1,#2)#3{\@xarg #1\relax \@yarg #2\relax
\@linelen=#3\unitlength
\ifnum\@xarg =0 \@vvector
  \else \ifnum\@yarg =0 \@hvector \else \@svector\fi
\fi}

\def\@hvector{\@hline\hbox to 0pt{\@linefnt
\ifnum \@xarg <0 \@getlarrow(1,0)\hss\else
    \hss\@getrarrow(1,0)\fi}}

\def\@vvector{\ifnum \@yarg <0 \@downvector \else \@upvector \fi}

\def\@svector{\@sline
\@tempcnta\@yarg \ifnum\@tempcnta <0 \@tempcnta=-\@tempcnta\fi
\ifnum\@tempcnta <5
  \hskip -\wd\@linechar
  \@upordown\@clnht \hbox{\@linefnt  \if@negarg
  \@getlarrow(\@xarg,\@yyarg) \else \@getrarrow(\@xarg,\@yyarg) \fi}%
\else\@badlinearg\fi}

\def\@getlarrow(#1,#2){\ifnum #2 =\z@ \@tempcnta='33\else
\@tempcnta=#1\relax\multiply\@tempcnta \sixt@@n \advance\@tempcnta
-9 \@tempcntb=#2\relax\multiply\@tempcntb \tw@
\ifnum \@tempcntb >0 \advance\@tempcnta \@tempcntb\relax
\else\advance\@tempcnta -\@tempcntb\advance\@tempcnta 64
\fi\fi\char\@tempcnta}

\def\@getrarrow(#1,#2){\@tempcntb=#2\relax
\ifnum\@tempcntb < 0 \@tempcntb=-\@tempcntb\relax\fi
\ifcase \@tempcntb\relax \@tempcnta='55 \or
\ifnum #1<3 \@tempcnta=#1\relax\multiply\@tempcnta
24 \advance\@tempcnta -6 \else \ifnum #1=3 \@tempcnta=49
\else\@tempcnta=58 \fi\fi\or
\ifnum #1<3 \@tempcnta=#1\relax\multiply\@tempcnta
24 \advance\@tempcnta -3 \else \@tempcnta=51\fi\or
\@tempcnta=#1\relax\multiply\@tempcnta
\sixt@@n \advance\@tempcnta -\tw@ \else
\@tempcnta=#1\relax\multiply\@tempcnta
\sixt@@n \advance\@tempcnta 7 \fi\ifnum #2<0 \advance\@tempcnta 64 \fi
\char\@tempcnta}

\def\@vline{\ifnum \@yarg <0 \@downline \else \@upline\fi}

\def\@upline{\hbox to \z@{\hskip -\@halfwidth \vrule
  width \@wholewidth height \@linelen depth \z@\hss}}

\def\@downline{\hbox to \z@{\hskip -\@halfwidth \vrule
  width \@wholewidth height \z@ depth \@linelen \hss}}

\def\@upvector{\@upline\setbox\@tempboxa\hbox{\@linefnt\char'66}\raise
     \@linelen \hbox to\z@{\lower \ht\@tempboxa\box\@tempboxa\hss}}

\def\@downvector{\@downline\lower \@linelen
      \hbox to \z@{\@linefnt\char'77\hss}}

\thinlines

\newcount\@xarg
\newcount\@yarg
\newcount\@yyarg
\newcount\@multicnt
\newdimen\@xdim
\newdimen\@ydim
\newbox\@linechar
\newdimen\@linelen
\newdimen\@clnwd
\newdimen\@clnht
\newdimen\@dashdim
\newbox\@dashbox
\newcount\@dashcnt
 \catcode`@=12

\newbox\tbox
\newbox\tboxa

\def\leftzer#1{\setbox\tbox=\hbox to 0pt{#1\hss}%
     \ht\tbox=0pt \dp\tbox=0pt \box\tbox}

\def\rightzer#1{\setbox\tbox=\hbox to 0pt{\hss #1}%
     \ht\tbox=0pt \dp\tbox=0pt \box\tbox}

\def\centerzer#1{\setbox\tbox=\hbox to 0pt{\hss #1\hss}%
     \ht\tbox=0pt \dp\tbox=0pt \box\tbox}

\def\image(#1,#2)#3{\vbox to #1{\offinterlineskip
    \vss #3 \vskip #2}}

\def\leftput(#1,#2)#3{\setbox\tboxa=\hbox{%
    \kern #1\unitlength
    \raise #2\unitlength\hbox{\leftzer{#3}}}%
    \ht\tboxa=0pt \wd\tboxa=0pt \dp\tboxa=0pt\box\tboxa}

\def\rightput(#1,#2)#3{\setbox\tboxa=\hbox{%
    \kern #1\unitlength
    \raise #2\unitlength\hbox{\rightzer{#3}}}%
    \ht\tboxa=0pt \wd\tboxa=0pt \dp\tboxa=0pt\box\tboxa}

\def\centerput(#1,#2)#3{\setbox\tboxa=\hbox{%
    \kern #1\unitlength
    \raise #2\unitlength\hbox{\centerzer{#3}}}%
    \ht\tboxa=0pt \wd\tboxa=0pt \dp\tboxa=0pt\box\tboxa}

\unitlength=1mm

\spacedmath{2pt}
\vsize = 25truecm
\hsize = 16truecm
\hoffset = -.15truecm
\voffset = -.5truecm
\baselineskip15pt
\overfullrule=0pt
\def\moinss{\mathrel{\hbox{\kern1pt\vrule height 2.3pt depth -1.6pt
width 4.2pt\kern1pt}}}
\def\indp{\par\hskip0.5cm}

\def\ron{{\scriptstyle\circ}}
\font\klm=cmsy5
\def\dag{\hbox{\klm\char121}}
\null\vskip0.5cm
\centerline{\bf Conformal blocks, fusion rules and the Verlinde formula}
\smallskip \centerline{Arnaud {\pc BEAUVILLE} \footnote{\parindent
0.6cm (*)}{\vtop{\eightpoint\baselineskip12pt\hsize15.2truecm\noindent
 Partially supported by the
European Science project ``Geometry of Algebraic Varieties", Contract
 SCI-0398-C(A).}}}\vskip0.7cm

\hfill{\it Dedicated to  F. Hirzebruch} \vskip0.7cm

{\it Introduction}
\smallskip

\ind The Verlinde formula  computes the dimension of certain vector spaces, the
{\it spaces of
conformal blocks},  which are the basic objects of a particular kind of quantum
field theories, the
so-called Rational Conformal Field Theories (RCFT). These spaces appear as
spaces of global multiform
sections of some flat vector bundles on the moduli space of curves with marked
points, so that their
dimension is simply the rank of the corresponding vector bundles.
  The
computation relies  on the behaviour of these bundles under degeneration of the
Riemann surface,
often referred to as the {\it factorization rules.} Verlinde's derivation from
the formula [V]
rested on a conjecture which does not seem to be proved yet in this very
general framework. \ind The
Verlinde formula has
 attracted a great deal of attention from the mathematicians when it was
realized that for
some particular RCFT's  associated to a compact Lie group $G$ (the WZW-models),
the spaces of conformal
blocks had a nice interpretation as spaces of {\it generalized theta
functions}, that is sections of
a determinant bundle (or its tensor powers) over the moduli space of $G$\tx
bundles on a Riemann
surface. This interpretation has been worked out  rigorously  for $SU(n)$ in
[B-L], and for the
general case in [F], while the factorization rules for these models
have been established  in [T-U-Y] and also in [F]. However there seems to be
some confusion
among mathematicians as to whether this work implies  the explicit Verlinde
formula  for the
spaces of generalized theta functions or not -- perhaps because of a few
misprints and inadequate
references in some of the above quoted papers.
 \ind Thus the aim of this paper is
to explain how the Verlinde formula for the WZW-models (hence for the space of
generalized
theta functions) can be derived from the factorization rules, at least in the
$SU(n)$ case.  As the
title indicates, the paper has three parts.  In the first one, which is
probably the most involved
technically, we fix a simple Lie algebra ${\goth g}$;  following [T-U-Y] we
associate a vector space
$V_C(\vec P,\vec \lambda)$ to a Riemann surface $C$ and a finite number of
points of $C$, to each of
which is attached a representation of ${\goth g}$. The main novelty here is a
more concrete
interpretation of this space (prop. 2.3) which gives a simple expression in the
case $C={\bf P}^1$ --
an essential ingredient of the Verlinde formula.  In the second part we develop
the formalism of the
{\it fusion rings},  an elegant
 way of encoding the factorization rules; this gives an explicit formula for
the dimension of
$V_C(\vec P,\vec \lambda)$ in terms of the characters of the fusion ring.
 In the third part we apply this formalism to the special case considered in
part I; this leads
to the fusion ring ${\cal R}_\ell({\goth g})$ of representations of level $\le
\ell$. We show
following [F] how one can determine the characters of ${\cal R}_\ell({\goth
g})$ when ${\goth g}$ is
${\goth sl}(n,{\bf C})$ or ${\goth sp}(n,{\bf C})$ (Faltings handles all the
classical algebras and
$G_2$, but there  seems to be no proof for the other exceptional algebras).
Putting things together
we obtain in these cases the Verlinde formula for the dimension of  $V_C(\vec
P,\vec \lambda)$.
\ind   I have tried  to make the paper as self-contained as possible, and in
particular not to assume
that the reader is an expert in Kac-Moody algebras; however some familiarity
with classical Lie
theory will certainly help. I would like to mention the
preprint [S] which contains (among other things) results  related to our Parts
II and III -- though
with a slightly different point of view. \vtop{\eightpoint\baselineskip12pt\ind
 I would like to thank
Y. Laszlo, O. Mathieu and C. Sorger for useful discussions.} \vskip1.2cm
\centerline{\bf Part I: the spaces $V^{}_C({\vec P},\vec\lambda)$} \smallskip
1. {\it Affine Lie
algebras }
 \ind (1.1) Throughout this paper we fix  a simple complex Lie algebra ${\goth
g}$, and  a Cartan
subalgebra  ${\goth h}\i{\goth g}$. I refer e.g. to [Bo] for the definition of
the root system
$R({\goth g,h})\i{\goth h}^*$, and of the coroot $H_\alpha\in{\goth
h}$ associated to a root $\alpha$. We have a decomposition ${\goth g}={\goth
h}\oplus\sum_{\alpha\in R({\goth g,h})}{\goth g}^\alpha$.   We
also fix a basis $(\alpha_1,\ldots,\alpha_r)$ of the root system, which
provides us with a
partition of the roots into positive and negative ones.
\ind  The {\it weight lattice} $P\i {\goth h}^*$
is the group of linear forms $\lambda\in{\goth h}^*$ such that
$\lambda(H_\alpha)\in{\bf Z}$ for all
roots $\alpha$. A weight $\lambda$ is {\it dominant} if
$\lambda(H_\alpha)\ge 0$ for all positive roots $\alpha$; we denote by $P_+$
the set of dominant
weights.  To each dominant weight $\lambda$ is associated a simple ${\goth
g}$\tx module $V_\lambda$,
unique up to isomorphism, containing a {\it highest weight vector} $v_\lambda$
with weight $\lambda$
(this means that $v_\lambda$ is annihilated by ${\goth g}^\alpha$ for
$\alpha>0$ and   that
$H\,v_\lambda=\lambda(H)v_\lambda$ for all $H$ in ${\goth h}$). The map
$\lambda\mapsto[V_\lambda]$ is a bijection of $P_+$ onto the set of isomorphism
classes of
finite-dimensional simple ${\goth g}$\tx modules.

\ind  (1.2) The {\it normalized Killing form} $(\ |\ )$ on ${\goth g}$ is
the unique
 ${\goth g}$\tx invariant nondegenerate symmetric form on ${\goth g}$
satisfying
$(H_\beta\,|\,H_\beta)=2$ for every long root ${\beta}$. We'll denote by the
same symbol the
non-degenerate form induced on ${\goth h}$ and the inverse form on ${\goth
h}^*$. We will  use these
normalized forms throughout the paper.
\ind (1.3) Let $\theta$ be the highest root of $R({\goth g,h})$,
and $H_\theta$ the corresponding coroot. Following [Bo] we choose elements
$X_\theta$ in ${\goth
g}^\theta$ and $X_{-\theta}$ in ${\goth g^{-\theta}}$ satisfying
$$[H_\theta,X_{\theta}]=2X_{\theta}\quad,\quad
[H_\theta,X_{-\theta}]=-2X_{-\theta}\quad,\quad
[X_\theta,X_{-\theta}]=-H_\theta\ .$$
These elements span a Lie subalgebra ${\goth s}$ of ${\goth g}$, isomorphic to
${\goth sl}_2$, which
will play an important role in this paper.\smallskip
 \ind   (1.4) The affine Lie algebra $\widehat {\goth g}$
associated to ${\goth g}$ is   a central extension of ${\goth g}\otimes {\bf
C}((z))$ by ${\bf C}$:
$$\widehat{\goth g}= \bigl({\goth g}\otimes {\bf C}((z))\bigr)\oplus {\bf C}c\
,$$ the bracket of two
elements of  ${\goth g}\otimes {\bf C}((z))$ being given by  $$[X\otimes f,
Y\otimes g]=[X,Y]\otimes
fg\  +\  c\cdot (X\,|\,Y)\Res(g\,df)\ .$$  We denote by $\widehat {\goth g}_+$
and
$\widehat {\goth g}_-$ the  subspaces  ${\goth g}\otimes z{\bf C}[[z]]$ and
${\goth g}\otimes
z^{-1}{\bf C}[z^{-1}]$ of  $\widehat {\goth g}$, so that we have a
decomposition
$$\widehat {\goth g}=\widehat {\goth g}_-\oplus {\goth g}\oplus {\bf C}c\oplus
\widehat {\goth
g}_+\ .$$
 By the formula for the Lie bracket, each summand is  actually a Lie subalgebra
of $\widehat {\goth
g}$.

 \ind  (1.5) We fix an integer $\ell> 0$ (the
level); we are interested in the irreducible representations of $\widehat
{\goth g}$ which are {\it of
level} $\ell$, i.e. such that the central element $c$ of $\widehat {\goth g}$
acts as multiplication
by $\ell$. Let $P_\ell$ be the set of dominant weights $\lambda$  of ${\goth
g}$
such that $\lambda(H_\theta)\le \ell$. The fundamental result of the
representation theory of
$\widehat {\goth g}$ (see e.g. [K]) asserts that the reasonable representations
 of level $\ell$ are
classified by $P_\ell$. More precisely, for each $\lambda\in P_\ell$, there
exists a simple $\widehat
{\goth g}$\tx module  ${\cal H}_\lambda$ of level $\ell$, characterized up to
isomorphism by the
following property:\par
\centerline{\it The subspace of ${\cal H}_\lambda$ annihilated by $\widehat
{\goth g}_+$ is
isomorphic as a ${\goth g}$\tx module to $V_\lambda$.}

 In the sequel we will identify  $V_\lambda$ to the subspace of ${\cal
H}_\lambda$
annihilated by $\widehat {\goth g}_+$.\smallskip

\ind (1.6) We will need a few more technical details about the $\widehat
{\goth g}$\tx module  ${\cal H}_\lambda$. Let us first recall its construction.
Let ${\goth p}$ be
the Lie subalgebra ${\goth g}\oplus {\bf C}c\oplus \widehat {\goth
g}_+$ of $\widehat{{\goth g}}$. We extend the representation of ${\goth g}$ on
$V_\lambda$ by
letting $\widehat {\goth
g}_+$ act trivially and $c$ as $\ell\, {\rm Id}_{V_\lambda}$; we denote by
${\cal V}_\lambda$  the  induced  $\widehat{{\goth g}}$\tx module
$U(\widehat{{\goth
g}})\otimes_{U({\goth p})}V_\lambda$. It contains a unique
maximal $\widehat{{\goth g}}$\tx submodule ${\cal Z}_\lambda$; then  ${\cal
H}_\lambda$ is the
quotient ${\cal V_\lambda/{\cal Z}_\lambda}$.
\ind Since $U(\widehat{\goth g})$ is
isomorphic  as a $U(\widehat{{\goth g}}_-)$\tx module to $U(\widehat{\goth
g}_-)\otimes_{\bf C}U({\goth p})$, we see that
 {\it the natural map $U(\widehat{{\goth
g}}_-)\otimes_{\bf C}V_\lambda\longrightarrow {\cal V}_\lambda$ is an
isomorphism of
$\widehat{{\goth g}}_-$\tx modules.}

\ind Let us identify $V_\lambda$ to the submodule $1\otimes V_\lambda$ of
${\cal V}_\lambda$.
With the notation of (1.3), the  submodule ${\cal Z}_\lambda$ is generated by
the element
$(X_\theta\otimes z^{-1})^{\ell -\lambda(H_\theta)+1}\,v_\lambda$ ({\it
cf.}~[K], exerc. 12.12);
this element is annihilated by $\widehat{{\goth g}}_+$ (see remark (3.6)
below).
 \smallskip
\ind (1.7) An important observation (which plays a crucial role in conformal
field theory) is that
the representation theory of $\widehat{{\goth g}}$ is essentially independent
of the choice of the
local coordinate $z$. Let $u=u(z)$ be an element of ${\bf C}[[z]]$ with
$u(0)=0$, $u'(0)\not=0$. The
automorphism $f\mapsto f\rond u$ of ${\bf C}((z))$ induces an automorphism of
${\goth g}\otimes {\bf
C}((z))$, which extends to an automorphism  $\gamma_u$ of $\widehat {\goth g}$
(given by
$\gamma_u(X\otimes f)=X\otimes f\ron u$). Let $\lambda\in P_\ell$; since
$\gamma_u$ preserves
$\widehat {\goth g}_+$ and is the identity on ${\goth g}$, the representation
$\pi_\lambda\rond
\gamma_u$ is  irreducible, and the subspace annihilated by  $\widehat {\goth
g}_+$ is exactly
$V_\lambda$. Therefore the representation $\pi_\lambda\rond \gamma_u$ is
isomorphic to $\pi_\lambda$.
In other words, there is a canonical linear automorphism $\Gamma_u$ of ${\cal
H}_\lambda$ such that
$\Gamma_u\bigl((X\otimes f)v\bigr)=(X\otimes f\ron u)\ \Gamma_u(v)$ for $v\in
{\cal H}_\lambda$,
$X\otimes f\in\widehat {\goth g}$ and $\Gamma_u(v)=v$ for $v\in V_\lambda$.
\smallskip
 \ind (1.8) Let ${\goth
a}$ be a Lie algebra, $V$ a ${\goth a}$\tx module. The {\it space of
coinvariants} of $V$, denoted by
$[V]_{\goth a}$, is the largest  quotient of $V$ on which ${\goth a}$ acts
trivially, that is the
quotient of $V$ by the subspace spanned by the vectors $Xv$ for $X\in {\goth
a}$, $v\in V$. This is
also $V/U^+({\goth a})V$, where $U^+({\goth a})$ is the augmentation ideal of
$U({\goth a})$.
\ind  Let  $V$ and $W$ two ${\goth a}$\tx modules. Using  the canonical
anti-involution $\sigma$ of $U({\goth a})$ (characterized by $\sigma(X)=-X$ for
any $X$ in ${\goth
a}$) we can consider $V$ as a right $U({\goth a})$\tx module.  {\it Then the
space of coinvariants
$[V\otimes W]_{{\goth a}}$ is the tensor product $V\otimes_{U({\goth a})}W$}:
they are both  equal to
the quotient of $V\otimes W$ by the subspace spanned by the elements $Xv\otimes
w+v\otimes Xw$
($X\in{\goth a}$, $v\in V$, $w\in W$).

\vskip0.7cm
2. {\it The spaces $V^{}_C({\vec P},\vec\lambda)$}\smallskip

\ind (2.1) Let $C$ be a smooth, connected, projective curve  over ${\bf C}$.
For each affine
open set $U\i X$, we denote by ${\cal O}(U)$ the ring of algebraic functions on
$U$, and by ${\goth
g}(U)$ the Lie algebra ${\goth g}\otimes {\cal O}(U)$.

\ind We want to associate  a vector space to the data of $C$, a finite subset
$\vec P=\{P_1,\ldots,P_p\}$ of $C$, and an element $\lambda_i$ of $P_\ell$
attached to each $P_i$.
In order to do this we consider the $\widehat {\goth g}$\tx module ${\cal
H}_{\vec\lambda}:={\cal H}_{\lambda_1}\otimes\ldots\otimes{\cal
H}_{\lambda_p}$.
We choose
 a local coordinate $z_i$ at each $P_i$, and denote by $f^{}_{P_i}$ the Laurent
series at $P_i$ of an
element $f\in{\cal O}(C\moins\vec P)$. This  defines  for each $i$ a ring
homomorphism ${\cal
O}(C\moins\vec P)\longrightarrow {\bf C}((z))$, hence a Lie algebra
homomorphism
${\goth g}(C\moins\vec P)\longrightarrow {\goth g}\otimes {\bf C}((z))$. We
define an action of
${\goth g}(C\moins\vec P)$ on ${\cal H}_{\vec\lambda}$ by the formula
 $$(X\otimes f)\cdot
(v_1\otimes \ldots\otimes v_p)=\sum_iv_1\otimes \ldots\otimes(X\otimes
f^{}_{P_i})v_i\otimes
\ldots\otimes v_p\leqno(2.2)$$
 (that this is indeed a Lie algebra action follows from the residue formula,
which gives
$\sum_i\Res^{}_{P_i}f^{}_{P_i}dg^{}_{P_i}=0$).  Using the notation of (1.8), we
put
 $$V^{}_C(\vec
P,\vec\lambda)=[{\cal H}_{\vec\lambda}]_{{\goth g}(C\moinss\vec
P)}\qquad,\qquad V_C^{\dag}(\vec
P,\vec\lambda)=\Hom_{{\goth g}(C\moinss\vec P)}({\cal H}_{\vec\lambda},{\bf
C})\ ,$$ where ${\bf
C}$ is considered as a trivial ${\goth g}(C\moins\vec P)$\tx module. Of course
$V_C^{\dag}(\vec
P,\vec\lambda)$ is the dual of $V^{}_C(\vec P,\vec\lambda)$. By (1.7) these
spaces do not depend -- up
to a canonical isomorphism -- on the choice of the local coordinates
$z_1,\ldots,z_p$. On the other
hand it is important to keep in mind that they depend on the Lie algebra
${\goth g}$ and the integer
$\ell$, though neither of these  appear in the notation.
\ind Though this will play no role in the sequel, I would like to mention that
these spaces have a
natural interpretation in the framework of algebraic geometry. Let me restrict
for simplicity to
the case ${\goth g}={\goth sl}_r({\bf C})$. Then the space
$V^{\dag}_C(\emptyset)$ is canonically
isomorphic to $H^0({\cal SU}_C(r),{\cal L}^\ell )$, where ${\cal SU}_C(r)$ is
the moduli space of
semi-stable vector bundles on $C$ with trivial determinant on $C$ and ${\cal
L}$ the
determinant line bundle (see [B-L], and [F] for the case of an arbitrary simple
Lie algebra). A
similar interpretation for $V^{\dag}_C(\vec P,\vec\lambda)$ has been worked out
by C.~Pauly in terms
of moduli spaces of parabolic vector bundles.  \medskip {\pc PROPOSITION}
2.3$.-$ {\it Let} ${\vec
P}=\{P_1,\ldots,P_p\}$, ${\vec Q}=\{Q_1,\ldots,Q_q\}$ {\it be two finite
non\-empty subsets of $C$,
without common point; let $\lambda_1,\ldots,\lambda_p${\rm ;} $
\mu_1,\ldots,\mu_q$ be elements of
$P_\ell$. We let ${\goth g}(C\moins {\vec P})$ act on $V_{\mu_j}$ through the
evaluation map
$X\otimes f\mapsto f(Q_j)X$. The inclusions $V_{\mu_j}\mono {\cal H}_{\mu_j}$
induce an isomorphism}
$$[{\cal H}_{\vec\lambda}\otimes V_{\vec \mu}]_{{\goth g}(C\moinss {\vec
P})}\iso[{\cal
H}_{\vec\lambda}\otimes{\cal H}_{\vec\mu}]_{{\goth g}(C\moinss{\vec
P}\moinss{\vec Q})}=V_C^{}(\vec
P\cup\vec Q,(\vec\lambda,\vec\mu))\ .$$ \smallskip \ind The case $\vec
Q=\{Q\}$, $\mu=0$ gives:
\smallskip {\pc COROLLARY} 2.4$.-$ {\it Let $Q\in C\moins {\vec P}$. There is a
canonical
isomorphism} $$V^{}_C({\vec P},\vec\lambda)\iso V^{}_C(\vec P\cup \vec
Q,(\vec\lambda,0))\
\quad\carre$$ \ind This is the ``propagation of vacua", {\it cf.}
 [T-U-Y], prop. 2.2.3.
\smallskip
{\pc COROLLARY} 2.5$.-$ {\it Let $Q\in C\moins {\vec P}$. There is a canonical
isomorphism}
$$V^{}_C({\vec P},\vec\lambda)\iso [{\cal H}_0\otimes V_{\vec\lambda}]_{{\goth
g}(C\moinss Q)}\ .$$
\ind Apply cor. 2.4, then the proposition inverting the role of $\vec P$ and
$\vec Q$. \cqfd
\smallskip
\ind (2.6) If $\vec\lambda=(0,\ldots,0)$, cor. 2 shows that $V^{}_C({\vec
P},\vec\lambda)$
is canonically isomorphic to $[{\cal H}_0]_{{\goth g}(C\moinss Q)}$, and in
particular independent of
$\vec P$. It follows that the space $[{\cal H}_0]_{{\goth g}(C\moinss Q)}$ is
independent of $Q$ up
to a canonical isomorphism;  we'll denote it by  $V^{}_C(\emptyset)$. Note that
with this convention
cor. 1 still holds in the case $\vec P=\emptyset$. \bigskip

 \ind I believe that the expression for  $V^{}_C({\vec P},\vec\lambda)$ given
by cor. 2.5 is more
flexible than the original definition. For instance we are going to use it
below to get a more
explicit expression in the case $C={\bf P}^1$. Also an easy proof of the
``factorization rules"
([T-U-Y], prop. 2.2.6) can be given in this set-up.
\smallskip
\ind (2.7) Let me finish with an easy result which we will need later on. For
each $\lambda\in
P_+$, the dual $V_\lambda^*$ is a simple  ${\goth g}$\tx module; let us denote
by $\lambda^*$ its
highest weight. The map $\lambda\mapsto \lambda^*$ is an involution of $P_+$,
which is actually the
restriction of a ${\bf Z}$\tx linear involution of $P$ (the experts have
already recognized the
automorphism $-w_0$, where $w_0$ is the  element of biggest length in the Weyl
group). This
involution also induces an invoIution of the root system which preserves the
root system, its
basis, and therefore the longest root $\theta$.  An important consequence is
that $P_\ell $ {\it is
preserved by the involution} $\lambda\mapsto \lambda^*$. \smallskip

\smallskip {\pc PROPOSITION} 2.8$.-$ {\it Put
$\vec\lambda^*=(\lambda_1^*,\ldots,\lambda_p^*)$. There is a natural
isomorphism
$$V^{}_C({\vec P},\vec\lambda)\iso V^{}_C({\vec P},\vec\lambda^*)\
$$}\vskip-4pt
\ind (This isomorphism is canonical once certain choices (a ``Chevalley basis")
have been done for the
Lie algebra ${\goth g}$.)\ind There exists an automorphism $\sigma$ of ${\goth
g}$ such that for each
finite-dimensional representation $\rho:{\goth g}\longrightarrow \End(V)$,
$\rho\rond \sigma$ is
isomorphic to the dual representation ([Bo], ch. VIII, \S 7, \no 6, remarque
1). The automorphism
$\sigma$ extends to an automorphism $\hat\sigma$ of $\widehat{{\goth g}}$,
which preserves the
decomposition $\widehat {\goth g}=\widehat {\goth g}_-\oplus {\goth g}\oplus
{\bf C}c\oplus \widehat {\goth
g}_+$.
 \ind Let $\lambda\in P_\ell $, and let $\pi_\lambda:\widehat{{\goth
g}}\longrightarrow \End({\cal
H}_\lambda)$ be the corresponding representation. The representation
$\pi_\lambda\rond \hat\sigma$
is simple, the subspace of ${\cal H}_\lambda$ annihilated by $\widehat{{\goth
g}}_+$ is
$V_\lambda$, on which ${\goth g}$ acts by the representation $\rho_\lambda\rond
\sigma$; therefore
$\pi_\lambda\rond \hat\sigma$ is isomorphic to $\pi_{\lambda^*}$. In other
words, there exists for
each $\lambda\in P_\ell $ a ${\bf C}$\tx linear isomorphism $t_\lambda:{\cal
H}_\lambda\longrightarrow {\cal H}_{\lambda^*}$ such that
$t_\lambda(Xv)=\hat\sigma(X)v$ for
$X\in\widehat{{\goth g}}$, $v\in{\cal H}_\lambda$.
\ind Now let
 $t_{\vec\lambda}:{\cal H}_{\vec\lambda}\longrightarrow {\cal
H}_{\vec\lambda^*}$ be the ${\bf C}$\tx
linear isomorphism $t_{\lambda_1}\otimes\ldots\otimes
t_{\lambda_p}$. It follows  from (2.2) and the above formula  that
$t_{\vec\lambda}\bigl((X\otimes f)\,v\bigr)=(\sigma(X)\otimes f)\,
t_{\vec\lambda}(v)$ for
$X\in{\goth g}$, $f\in{\cal O}(X\moins \vec P)$, $v\in {\cal H}_{\vec\lambda}$.
Therefore
$t_{\vec\lambda}$ induces an isomorphism of
$V^{}_C({\vec P},\vec\lambda)$ onto $V^{}_C({\vec P},\vec\lambda^*)$. \cqfd

\vskip0.7cm 3. {\it Proof of Proposition} 2.3 \smallskip \ind Put $Q=Q_q$,
$\mu=\mu_q$, $U
=C\moins\vec P$, and
 ${\cal H}={\cal H}_{\vec\lambda}\otimes V_{\mu_1}\otimes\ldots\otimes
V_{\mu_{q-1}}$. Reasoning by induction on $q$ it will be enough to prove that
the inclusion
$V_\mu\mono {\cal H}_\mu$ induces an isomorphism $$[{\cal H}\otimes
V_{\mu}]_{{\goth
g}(U )}\iso[{\cal H}\otimes{\cal H}_{\mu}]_{{\goth g}(U \moinss Q)}\ .$$
\ind (3.1) Let me first explain the action of ${\goth g}(U \moins Q)$ on ${\cal
H}\otimes{\cal
H}_\mu$. We choose a local coordinate $z$ at $Q$. As
before, the map $X\otimes f\mapsto X\otimes f^{}_Q$ defines a Lie algebra
homomorphism
$\varepsilon:{\goth g}(U \moins Q)\longrightarrow {\goth g}\otimes{\bf
C}((z))$. Let us denote by
$\widehat {\goth g}(U \moins Q)$ the pull-back by $\varepsilon$ of the
extension $\widehat {\goth
g}\longrightarrow {\goth g}\otimes {\bf C}((z))$; in other words, $\widehat
{\goth g}(U \moins Q)$ is
the space ${\goth g}(U \moins Q)\oplus{\bf C}c$, the bracket of two elements
$X\otimes f$, $Y\otimes g$
being given by
$$[X\otimes f, Y\otimes g]=[X,Y]\otimes fg\
+\  c\cdot (X\,|\,Y)\Res^{}_Q(g\,df)\ .$$
Applying again the Residue formula we see that the action of ${\goth g}(U
\moins Q)$ on ${\cal
H}_{\vec\lambda}$ given by formula (2.2) extends to an action of $\widehat
{\goth g}(U \moins Q)$,
which is {\it of level} $-\ell$ in the sense that the central element $c$ acts
as multiplication by
$-\ell$. On the other hand  $\varepsilon$ extends by construction to a
homomorphism $\widehat {\goth
g}(U \moins Q)\longrightarrow \widehat {\goth g}$ through which $\widehat
{\goth g}(U \moins Q)$ acts on
${\cal H}_\mu$ with level $\ell$, hence the action on ${\cal H}\otimes{\cal
H}_\mu$ is of level $0$
and therefore factors through ${\goth g}(U \moins Q)$.

\ind Besides the fact that it is of level $-\ell $, the only property we will
use of the action of
$\widehat{\goth g}(U \moins Q)$ on ${\cal H}$ is the following:

\centerline{(*)  {\it The endomorphism $X_{-\theta}\otimes f$ of ${\cal H}$ is
locally nilpotent for
all} $f\in {\cal O}(U)$. }
(This is because $X_{-\theta}$ is a nilpotent element of ${\goth g}$, while
every element of
$\widehat{{\goth g}}_+$ is locally nilpotent in the integrable modules ${\cal
H}_{\lambda_i}$.)
 \bigskip

\ind (3.2) We first check that the map ${\cal
H}\otimes V_\mu\mono {\cal H}\otimes {\cal H}_\mu $ is equivariant with respect
to ${\goth
g}(U )$. This amounts to prove that the inclusion $V_\mu \mono {\cal H}_\mu $
is equivariant.
But $V_\mu$ is the subspace of ${\cal H}_\mu$ annihilated by $\widehat{\goth
g}_+$ (1.5), so an
element $X\otimes f$ of ${\goth g}(U )$ acts on ${\cal H}_\mu$ as the element
$f(Q)X$ of ${\goth
g}$, hence our assertion. Therefore the inclusion induces a linear map
$$i: [{\cal H}\otimes V_\mu]_{{\goth g}(U )}\longrightarrow  [{\cal
H}\otimes {\cal H}_\mu]_{{\goth g}(U \moinss Q)}\ .$$

\ind (3.3) We prove that the statement is true when we replace the simple
module ${\cal
H}_\mu$ by the  module ${\cal V}_\mu$ (1.6).
 Let us observe first that by (1.7), the statement is independent of the choice
of the local
coordinate $z$ at $Q$. We choose $z$ so that  $z^{-1}\in {\cal O}(U \moins Q)$
(this is possible as
soon as $\vec P \not=\emptyset$). From the decomposition $\displaystyle {\cal
O}(U \moins Q)={\cal
O}(U )\oplus \sum_{n\ge 1}{\bf C}z^{-n}$ we get
$${\goth g}(U \moins Q)={\goth g}(U )\oplus\widehat {\goth g}_-$$
where we have identified the Lie algebra  $\sum_{n\ge
1}{\goth g}\ z^{-n}$ with its image $\widehat {\goth g}_-$ in $\widehat {\goth
g}$. Note that both
summands can be viewed as Lie subalgebras of $\widehat {\goth g}(U \moins Q)$.
\ind Let us consider
first the coinvariants under $\widehat {\goth g}_-$. By (1.8)
 $[{\cal H}\otimes {\cal V}_\mu]_{\widehat {\goth g}_-}$ can be identified with
${\cal H}\otimes_{U(\widehat {\goth g}_-)} {\cal V}_\mu$. Since the natural map
$\widehat {\goth g}_-\otimes_k V_\mu\longrightarrow {\cal V}_\mu$ is an
isomorphism
of $\widehat {\goth g}_-$\tx modules (1.6), we conclude that the inclusion
$V_\mu\mono{\cal V}_\mu$
induces an isomorphism ${\cal H}\otimes V_\mu\iso  [{\cal
H}\otimes {\cal V}_\mu]_{\widehat {\goth g}_-}$. Taking coinvariants under
${\goth g}(U )$ gives
the required isomorphism.
\smallskip
\ind (3.4)  Let ${\cal
Z}_\mu$ be the kernel of the canonical surjection ${\cal V}_\mu\rightarrow
{\cal
H}_\mu$; we have an exact sequence
$${\cal H}\otimes {\cal Z}_\mu\longrightarrow [{\cal H}\otimes{\cal
V}_\mu]_{{\goth g}(U \moinss Q)}\longrightarrow  [{\cal H}\otimes{\cal H}_\mu
]_{{\goth g}(U \moinss Q)}\rightarrow 0\ ,$$
so we want to prove that the image of ${\cal H}\otimes {\cal Z}_\mu$ in ${\cal
H}
\otimes_{U(\widehat {\goth g}(U \moinss Q))}{\cal V}_\mu$ is zero. As a
$U(\widehat {\goth g})$\tx
 module ${\cal Z}_\mu$ is generated by the vector $(X_\theta\otimes
z^{-1})^{k}\,v_\mu$, where $v_\mu$
is a highest weight vector and $k=\ell-\mu(H_\theta)+1$ (1.6); moreover this
vector is annihilated
by $\widehat{{\goth g}}_+$, so it generates ${\cal Z}_\mu$ as a
$U(\widehat{\goth g}_-\oplus{\goth
g})$\tx
 module. Since $\widehat {\goth g}_-\oplus{\goth
g}\i {\goth g}(U \moins Q)$, it is enough to prove that $h\otimes
(X_\theta\otimes
z^{-1})^{k}\,v_\mu=0$ in  ${\cal H} \otimes_{U({\goth g}(U \moinss Q))}{\cal
V}_\mu$ for each vector
$h\in{\cal H}$. Let $f$ be an element of ${\cal O}(U )$ such that $f^{}_Q\equiv
z$ (mod.$\,z^2$); put
$Y=X_{-\theta}\otimes f$. By property $(*)$ in (3.1)  there exists an integer
$N$ such that $Y^Nh=0$.
By  lemma (3.5) below $(X_\theta\otimes z^{-1})^{k}\,v_\mu$ can be written
$Y^Nw$ for some $w\in
{\cal V}_\mu$, so $h\otimes (X_\theta\otimes z^{-1})^{k}\,v_\mu$ is zero in
${\cal H}
\otimes_{U(\widehat {\goth g}(U \moinss Q))}{\cal V}_\mu$, which finishes the
proof.

\medskip
 {\it Lemma} 3.5$.-$ {\it Let
$f(z)\in {\bf C}[[z]]$ such that $f(0)=0$, $f'(0)=1$. Put $X=X_\theta\otimes
z^{-1}$,
$Y=X_{-\theta}\otimes f(z)$ in $\widehat {\goth g}$. Let $\mu\in P_\ell$, and
$p,q\in{\bf N}$ with
$p\ge \ell+1-\mu(H_\theta)$. There exists a nonzero rational number
$\alpha_{p,q}$ such that
$X^{p}v_\mu=\alpha_{p,q}Y^{q}X^{p+q}v_\mu$ in ${\cal V}_\mu$.}
\ind Let $H:=[Y,X]=\bigl(H_\theta\otimes z^{-1}f(z)\bigr)-c$; then
$[H,X]=2X_\theta\otimes z^{-2}f(z)$ commutes with $X$, so one has  in
$U(\widehat {\goth g})$
$$HX^m=X^mH+\sum_{a+b=m-1}X^a[H,X]X^b=X^mH+mX^{m-1}[H,X]$$ Since $v_\mu$ is
annihilated by $\widehat
{\goth g}_+$ and by $X_\theta$, one has $[H,X]v_\mu=2Xv_\mu$ and
$Hv_\mu=-kv_\mu$ with
$k=\ell-\mu(H_\theta)$, hence $HX^mv_\mu=(2m-k)X^mv_\mu$. Then
$$YX^{p+1}v_\mu=\sum_{n+m=p}X^nHX^mv_\mu=(p+1)(p-k)X^{p}v_\mu\ .$$ This proves
the lemma
in the case $q=1$; the general case follows at once by induction on $q$. \cqfd
\medskip
{\it Remark} 3.6$.-$ The same method gives the vanishing of
$YX^{\ell-\mu(H_\theta) +1}\,v_\mu$ for
any $Y\in \widehat{{\goth g}}_+$.

\vskip0.7cm

4. {\it The case $C={\bf P}^1$}
\smallskip
\ind We fix a coordinate $t$ on ${\bf P}^1$.
\medskip
{\pc PROPOSITION} 4.1$.-$ {\it Let $P_1,\ldots,P_p$ be distinct points of ${\bf
P}^1$, with coordinates $t_1,\ldots,t_p$, and let $\lambda_1,\ldots,\lambda_p$
be elements of
$P_\ell$. Let $T$ be the endomorphism of $V_{\vec\lambda}$ defined by
$$T(v_1\otimes\ldots\otimes v_p)=\sum_{i=1}^p t_i\, v_1\otimes\ldots\otimes
 X_\theta\, v_i \otimes\ldots\otimes v_p
\ .$$ The space
 $V_{{\bf P}^1}(\vec P,\vec\lambda)$ is canonically isomorphic to the
largest quotient of $V_{\vec\lambda}$ on which ${\goth g}$ and
$T^{\ell+1}$ act trivially. The space $V^{\dag}_{{\bf P}^1}(\vec
P,\vec\lambda)$
is isomorphic to the space of ${\goth g}$\tx invariant $p$\tx linear forms
$\varphi:
V_{\lambda_1}\times\ldots\times V_{\lambda_p}\longrightarrow {\bf C}$ such
that} $\varphi\rond
T^{\ell+1}=0$.
\ind We apply cor. 2.5 with $Q=\infty$ (so that the local coordinate $z$ at $Q$
is $t^{-1}$). This
gives an isomorphism of  $V_{{\bf P}^1}(\vec P,\vec\lambda)$ onto $[{\cal
H}_0\otimes
V_{\vec\lambda}]_{{\goth g} ({\bf A}^1)}$. Now ${\goth g}({\bf A}^1)$ is the
sum of ${\goth g}$ and
$\widehat {\goth g}_-$; it follows from (1.6) that the  $U({\goth g}({\bf
A}^1))$\tx module ${\cal
H}_0$ is generated by the highest weight vector $v_0$, with the relations
${\goth g}\,v_0=0$ and
$(X_\theta\otimes z^{-1})^{\ell+1}\,v_0=0$. Therefore the space $[{\cal
H}_0\otimes
V_{\vec\lambda}]_{{\goth g}({\bf A}^1)}\cong {\cal H}_0\otimes_{U({\goth
g}({\bf A}^1))
}V_{\vec\lambda}$ is canonically isomorphic to  $V_{\vec\lambda}/({\goth
g}V_{\vec\lambda}+\im
T^{\ell +1})$, where $T\ (=X_\theta\otimes t)$ is the endomorphism of
$V_{\vec\lambda}$ given by the
above formula.  The description of $V^{\dag}_{{\bf P}^1}(\vec P,\vec\lambda)$
follows by duality.
\cqfd

\bigskip
  \ind When $p=3$, one can describe the space $V_{{\bf
P}^1}(a,b,c;\lambda,\mu,\nu)$ (or its dual) in a
more concrete way. Let us first consider the case when ${\goth g}={\goth
sl}_2$. We denote by $E$ the
standard 2-dimensional representation of ${\goth g}$. We will identify $P_\ell$
 with the set  of integers $p$ with $0\le p\le \ell$ (by associating to such an
integer the representation $S^pE$). By prop. 2, $V^{\dag}_{{\bf
P}^1}(a,b,c;p,q,r)$ is the space
of linear forms $F\in\Hom_{\goth g}(S^pE\otimes S^qE\otimes S^rE,{\bf C})$ such
that $F\rond
T^{\ell+1}=0$. \medskip
{\it Lemma} 4.2$.-$ a) {\it  The space
$\Hom_{\goth g}(S^pE\otimes S^qE\otimes S^rE,{\bf C})$  is either $0$ or $1$\tx
dimensio\-nal.
It is nonzero if and only if $p+q+r$ is even, say $=2m$,
and $p,q,r$ are} $\le m$.
\ind  b)  {\it The subspace $V^{\dag}_{{\bf P}^1}(a,b,c;p,q,r)$ is nonzero if
and only if $p+q+r$ is
even and} $\le 2\ell$.
 \ind The first assertion is an immediate consequence of the
Clebsch-Gordan formula. When the space $\Hom_{\goth g}(S^pE\otimes S^qE\otimes
S^rE,{\bf C})$ is
nonzero, a generator  $G$ is obtained as follows: the dual of $S^pE\otimes
S^qE\otimes
S^rE$ is the space of polynomial maps $F(u,v,w)$ \ $(u,v,w\in E)$ which are
homogeneous of degree $p$
in $u$, $q$ in $v$ and $r$ in $w$; then the polynomial
$G(u,v,w)=\varphi(v,w)^{n-p}\varphi(w,u)^{n-q}\varphi(u,v)^{n-r}$, where
$\varphi$ is any non-zero
alternate form on $E$, is  clearly ${\goth g}$\tx invariant. It remains to make
explicit
the action of $T$ on the dual of $S^pE\otimes S^qE\otimes S^rE$. \ind This dual
can also be seen  as
the space $P_{p,q,r}$ of (non homogeneous) polynomials $P(x,y,z)$ of degree
$\le p$ in $x$, $\le q$
in $y$ and $\le r$ in $t$: the correspondence is obtained by choosing a basis
$(e_0,e_1)$ of $E$ and
putting $P(x,y,z)=$ $F(e_0+xe_1,e_0+ye_1,e_0+ze_1)$. In particular, the
polynomial corresponding to
$G$ is (up to a constant)  $Q(x,y,z)=(y-z)^{n-p}(z-x)^{n-q}(x-y)^{n-r}$.
\ind  Choose the basis so that
$X_\theta \,e_0=e_1$; then the action of $X_\theta\otimes 1\otimes 1$ (resp.
$1\otimes X_\theta\otimes
1$, resp. $1\otimes 1\otimes X_\theta$) on $P_{p,q,r}$ is the derivation with
respect to $x$ (resp.
$y$, resp. $z$). Therefore  $T$ acts as the operator $a{\partial\over \partial
x}+b{\partial\over
\partial y}+c{\partial\over \partial z}$; in other words, $T^m\cdot P$ is the
coefficient of
${h^m\over m!}$ in the expansion of $P(x+ah,y+bh,z+ch)$. Since
$Q(x+ah,y+bh,z+ch)$ is a polynomial
of degree $n$ in $h$ (because $a\not=b\not=c$), we obtain  b). \cqfd
\bigskip
\ind (4.3) In the general case, we consider
 the  Lie subalgebra  ${\goth s}\cong{\goth sl}_2$ of ${\goth g}$  with basis
$(X_\theta,X_{-\theta},H_{\theta})$ (1.3). Following the (unpleasant) practice
of the physicists,
we'll say that an irreducible representation of ${\goth sl}_2$ {\it has spin}
$i$ if it is isomorphic
to $S^{2i}E$; so the spin is a half-integer. Let $\lambda\in P_\ell$; as a
${\goth s}$-module, $V_\lambda$ breaks as a direct sum of isotypic components
$V_\lambda^{(i)}$ of
spin $i$, with $0\le i\le \ell /2$.
 \smallskip {\pc PROPOSITION} 4.3$.-$ a) {\it The space $V_{{\bf
P}^1}(a,b,c;\lambda,\mu,\nu)$ is canonically isomorphic to the quotient of
$[V_{\lambda}\otimes
V_\mu\otimes V_\nu]_{\goth g}$ by the image of the subspaces
$V_{\lambda}^{(p)}\otimes
V^{(q)}_\mu\otimes V^{(r)}_\nu$ for $p+q+r>\ell$.}
\ind b) {\it
The space $V^{\dag}_{{\bf P}^1}(a,b,c;\lambda,\mu,\nu)$ is  canonically
isomorphic to the space of
${\goth g}$\tx invariant linear forms $\varphi: V_{\lambda}\otimes V_\mu\otimes
V_\nu\longrightarrow
{\bf C}$ which vanish on the subspaces $V_{\lambda}^{(p)}\otimes
V^{(q)}_\mu\otimes V^{(r)}_\nu$
whenever $p+q+r>\ell$.}

\ind The two assertions are of course equivalent; let us prove b). By prop.
4.1, all we have to do is
to express the condition $\varphi\rond T^{\ell+1}$ for a ${\goth g}$\tx
invariant linear form
$\varphi: V_{\lambda}\otimes V_\mu\otimes V_\nu\longrightarrow {\bf C}$. Write
$V_\lambda=\build\oplus_{p=0}^{\ell/2} V_\lambda^{(p)}$, and similarly for
$V_\mu$ and $V_\nu$. The
subspaces $V_{\lambda}^{(p)}\otimes V^{(q)}_\mu\otimes V^{(r)}_\nu$ are stable
under ${\goth s}$ and
$T$, so we have to express that the restriction $\varphi_{pqr}$ of $\varphi$ to
any of these subspaces
vanishes on $\im T^{\ell+1}$. By the above lemma this is automatically
satisfied if $p+q+r\le \ell$,
while it  imposes $\varphi_{pqr}=0$ when $p+q+r> \ell$, hence the proposition.
\cqfd
\medskip
\ind Let me mention an easy consequence (which of course can also be proved
directly):
\smallskip
{\pc COROLLARY} 4.4$.-$ {\it One has}
$$\nospacedmath\displaylines{V_{{\bf P}^1}(P,\lambda)=0\quad{\rm
for}\quad\lambda\not=0\quad,\quad
V_{{\bf P}^1}(P,0)\cong V_{{\bf P}^1}(\emptyset)\cong {\bf C}\cr V_{{\bf
P}^1}(P,Q,\lambda,\mu)=0\quad{\rm for}\quad\mu\not=\lambda^*\quad,\quad V_{{\bf
P}^1}(P,Q,\lambda,\lambda^*)\cong{\bf C}\ .\quad\carre}$$
\medskip
\ind (4.5) Let me now recall one of the essential results of [T-U-Y], the {\it
factorization rules}
for the spaces $V^{}_C({\vec P},\vec\lambda)$.
In this paper we'll be only interested in the  dimension of these spaces, so I
will formulate
the factorization rules  in these terms. According to [T-U-Y] the dimension of
$V^{}_C({\vec P},\vec\lambda)$ depends only on
the genus $g$ of $C$ and of the set of weights
$\vec\lambda=(\lambda_1,\ldots,\lambda_p)$; let us
denote it by $N_g(\vec\lambda)$. One has
$$N_g(\vec\lambda)=\sum_{\nu\in P_\ell }N_{g-1}(\vec\lambda,\nu,\nu^*)\ ,$$
and, if
 $\vec\mu=(\mu_1,\ldots,\mu_q)$ is another set of weights and $h$, $k$
non-negative integers
such that $g=h+k$, $$N_g(\vec\lambda,\vec\mu)=\sum_{\nu\in P_\ell
}N_{h}(\vec\lambda,\nu)\,N_{k}(\vec\mu,\nu^*)\ .$$
\vskip1.2cm
\centerline{\bf Part II: Fusion rings}
\vglue 12pt plus 3pt minus 3pt
5. {\it Fusion rules and fusion rings}
\smallskip
\ind (5.1) Let $I$ be a finite set, with an involution $\lambda\mapsto
\lambda^*$. We'll denote by
${\bf N}^{(I)}$ the free commutative monoid generated by $I$, that is the set
of  sums
 $\displaystyle \sum_{\alpha\in I}n_\alpha\alpha$ with $n_\alpha\in {\bf N}$;
we shall always identify
$I$ to a subset of ${\bf N}^{(I)}$.
 The involution of $I$ extends by linearity to an involution $x\mapsto x^*$ of
${\bf N}^{(I)}$.
\smallskip
{\pc DEFINITION}$.-$  {\it A  fusion rule on $I$
is a map $N:{\bf N}^{(I)}\rightarrow {\bf Z}$ satisfying the following three
conditions:}
\indp (F 0) {\it One has $N(0)=1$, and  $N(\alpha)>0$ for some } $\alpha\in I$;
\indp (F 1) $N(x^*)=N(x)$ {\it  for every} $x\in {\bf N}^{(I)}$;
\indp (F 2) {\it For $x$, $y$ in ${\bf N}^{(I)}$, one has} $\displaystyle
N(x+y)=\sum_{\lambda\in
I}N(x+\lambda)\,N(y+\lambda^*)$. \medskip
\ind Let us call {\it kernel} of a fusion rule $N$ the set of elements $\alpha$
in $I$
such that $N(\alpha+x)=0$ for all $x\in {\bf N}^{(I)}$; one says that $N$ is
{\it non
degenerate} if its kernel is empty. Let $N$ be a  fusion rule  on $I$ with
kernel $K\i I$;  then $K$
is stable under $^*$, the restriction $N_0$ of $N$ to ${\bf N}^{(I\moinss K)}$
is a fusion rule on
$I\moins K$, and $N$ is simply the extension by $0$ of $N_0$ to  ${\bf
N}^{(I)}$. Therefore we can
restrict ourselves without loss of generality to the non-degenerate fusion
rules.

\medskip
{\it Examples} 5.2$.-$ {\it a}) Fix  a simple complex Lie algebra and a level
$\ell $. For
$\lambda_1,\ldots,\lambda_p$ in $P_\ell$, put using the notation of Part I
$$N(\sum
\lambda_i)=\dim V_{{\bf P}^1}(\vec P,\vec\lambda)\ ,$$ where $\vec P$ is an
arbitrary subset of ${\bf P}^1$ with $p$ elements. Then $N$ {\it is a
(non-degenerate) fusion rule on} $P_\ell$: the condition (F 0) and the
non-degeneracy condition
follows from cor. 4.4, (F 1) from prop. 2.8, and (F 2) is a particular case
$(h=k=0)$ of the
factorization rules (4.5).

\medskip
\ind {\it b}) Let $R$ be a commutative ring, endowed with an involutive ring
homomorphism $x\mapsto
x^*$ and a ${\bf Z}$\tx linear form $t:R\rightarrow {\bf Z}$; suppose
 that the bilinear form $(x,y)\mapsto t(xy^*)$ is symmetric and admits an
orthonormal basis
 $I$ (over ${\bf Z}$) containing $1$.  Define
a map $N: {\bf N}^{(I)}\longrightarrow {\bf Z}$ by the formula $N(\sum
n^{}_\alpha\alpha)=t(\prod \alpha^{n_\alpha})$. Then $N$ {\it is a
(non-degenerate) fusion rule on}
$I$. For the condition on $t$ implies in particular
$t(x^*)=t(x)$ and $t(1)=1$, hence
 (F 0) and (F 1). Since $I$ is an orthonormal basis, one has,
for $x$, $y$ in $R$,  $$t(xy)=\sum_{\lambda\in I}t(x\lambda)\,t(y\lambda^*)\
,\leqno(5.2)$$which implies  (F 2). Conversely:
\medskip
 {\pc PROPOSITION} 5.3$.-$ {\it Let $N:{\bf
N}^{(I)}\rightarrow {\bf Z}$ be a (non-degenerate) fusion rule on $I$. There
exists a ${\bf
Z}$\tx bilinear map ${\bf Z}^{(I)}\times{\bf Z}^{(I)}\longrightarrow {\bf
Z}^{(I)}$, which
turns ${\bf Z}^{(I)}$ into a  commutative ring, and a linear form $t$, uniquely
determined, such that
$$N(\sum
n^{}_\alpha\alpha)=t(\prod \alpha^{n_\alpha})$$for all elements $\sum
n_\alpha\alpha$ of ${\bf
N}^{(I)}$. One has
 $t(\alpha\beta^*)=\delta_{\alpha\beta}$ for $\alpha$, $\beta$ in $I$. }
\ind Let us apply (F 2) with $x=y=0$. Using
 (F 0) and (F 1) we get $\displaystyle \sum_{\lambda\in I}N(\lambda)^2=1$. This
means that there
exists an element $\varepsilon$ of $I$ such that
$$\varepsilon=\varepsilon^*\quad,\quad N(\varepsilon)=1\quad,\quad
N(\lambda)=0\quad {\rm
for}\quad\lambda\not=\varepsilon\ .$$ Then  (F 2) (with $y=0$) implies
$N(x+\varepsilon)=N(x)$ for
all $x\in {\bf N}^{(I)}$. \ind Now let us apply (F 2) with $x=\alpha$,
$y=\alpha^*$; we
obtain (using (F 1)) $$N(\alpha+\alpha^*)=\sum_{\lambda\in
I}N(\alpha+\lambda)^2\ \ge\
N(\alpha+\alpha^*)^2\ .$$If
 $N(\alpha+\lambda)=0$ for all $\lambda\in I$, one deduces from   (F 2)
$N(\alpha+x)=0$
for all $x\in {\bf N}^{(I)}$, which contradicts the non-degeneracy hypothesis.
Therefore the above
inequality implies $$ N(\alpha+\lambda)=0\quad{\rm
for}\quad \lambda\not=\alpha^* \quad,\quad N(\alpha+\alpha^*)=1\ .
\leqno(5.4)$$

\ind Let us define a multiplication law on ${\bf Z}^{(I)}$ by putting
$$\alpha\cdot\beta=\sum_{\lambda\in I}N(\alpha+\beta+\lambda^*)\,\lambda\
,\leqno(5.5)$$
and extending by bilinearity. This law is commutative; for $\alpha$, $\beta$,
$\gamma$ in $I$, one has acording to  (F 2)
$$(\alpha\cdot\beta)\cdot\gamma=\sum_{\lambda,\mu\in
I}N(\alpha+\beta+\lambda^*)\,N(\lambda+\gamma+\mu^*)\mu=\sum_{\mu\in
I}N(\alpha+\beta+\gamma+\mu^*)\,\mu=\alpha\cdot(\beta\cdot\gamma)\ ,$$so that
the  multiplication
is associative. One gets similarly, by induction on $s$,
$$\alpha_1\cdots\alpha_s=\sum_{\lambda\in
I}N(\alpha_1+\ldots+\alpha_s+\lambda^*)\
\lambda\leqno(5.6)$$
for $\alpha_1,\ldots,\alpha_s$ in $I$. Moreover one deduces from
(5.4)$$\varepsilon\cdot\alpha=\sum_{\lambda\in
I}N(\varepsilon+\alpha+\lambda^*)\,\lambda=\alpha\
$$Condition (F 1) implies that the involution $x\mapsto x^*$ is a ring
homomorphism.\ind  Let $t:{\bf
Z}^{(I)}\longrightarrow {\bf Z}$ be the  linear form $\sum
n_\alpha\alpha\mapsto n_\varepsilon$. One
gets from (3) $$t(\alpha_1\cdots\alpha_s)=N(\alpha_1+\ldots+\alpha_s)$$which is
the required formula
for $N$. Then (5.4) translates as $t(\alpha\beta^*)=\delta_{\alpha\beta}$.
\cqfd
\medskip
{\pc DEFINITION}$.-$ {\it The ring ${\bf Z}^{(I)}$ with the multiplication
given by $(5.5)$
is called the fusion ring associated to $N$.}
\ind We will denote it by ${\cal F}_N$ or simply ${\cal F}$.
\medskip
{\it Remark} 5.7$.-$ Most of the above still holds when one
replaces  (F 1) by the weaker condition

$(F\ 1')$: {\it One has $N(\alpha^*)=N(\alpha)$ and
$N(\alpha^*+\beta^*)=N(\alpha+\beta)$
for} $\alpha,\beta\in I$,

 the only difference being that the involution is not necessarily a ring
homomorphism -- in fact this
property is equivalent to (F 1).

\bigskip \ind (5.8) A consequence of prop. 5.3 is that the bilinear form
$(x,y)\mapsto t(xy)$ defines
an isomorphism of ${\cal F}$ onto the ${\cal F}$\tx module $\Hom_{\bf Z}({\cal
F},{\bf Z})$ (this
implies that ${\cal F}$ is a {\it Gorenstein} ${\bf Z}$\tx algebra). There is a
canonical element
$\Tr$ in $\Hom_{\bf Z}({\cal F},{\bf Z})$, the trace form: for each $x\in{\cal
F}$ the multiplication
by $x$ is an endomorphism $m_x$ of ${\cal F}$, and we put $\Tr(x):=\Tr(m_x)$.
An easy (and standard)
computation shows that the element of ${\cal F}$ which corresponds to $\Tr$
under the above
isomorphism is the {\it Casimir element} $\displaystyle
\omega:=\sum_{\lambda\in I}\lambda\lambda^*$;
in other words, one has $$\qquad\Tr(x)=t(\omega x)\quad\hbox{for every}\quad
x\in I\ .\leqno(5.8)$$
\medskip
 \ind  So far we have considered only fusion rules in genus $0$; this is no
restriction because the
general case reduces easily to the genus $0$ case. To be precise:
\smallskip
{\pc PROPOSITION} 5.9$.-$ {\it Suppose given a sequence of maps
$N_g:{\bf N}^{(I)}\longrightarrow {\bf Z}$ such that $N_0=N$ and
$$ N_g(x)=\sum_{\lambda\in
I}N_{g-1}(x+\lambda+\lambda^*)$$for $x$ in ${\bf N}^{(I)}$ and $g\ge 1$.
Then one has, for $g\ge 1$ and} $\alpha_1,\ldots,\alpha_p\in I$
$$N_g(\alpha_1+\ldots+\alpha_p)=t(\alpha_1\cdots\alpha_p\,\omega^g)=\Tr(\alpha_1\cdots\alpha_p\,
\omega^{g-1})\ .$$
\ind By induction on $g$ one gets
$$\eqalign{N_g(\alpha_1+\ldots+\alpha_p)&=\sum_{\lambda_1,\ldots,\lambda_g\in
I}N_0(\alpha_1+\ldots+\alpha_p+\lambda^{}_1+\lambda_1^*+\ldots+\lambda^{}_g+\lambda_g^*)\cr
&=\sum_{\lambda_1,\ldots,\lambda_g\in
I}t(\alpha_1\cdots\alpha_p\,\lambda^{}_1\lambda_1^*\cdots\lambda^{}_g\lambda_g^*)\cr
&=t(\alpha_1\cdots\alpha_p\,\omega^g)\ ;
}$$the last equality follows from (5.8).
\medskip
{\it Remark} 5.10$.-$ Using formula (5.2), it follows  that the sequence
$(N_g)$ satisfies
the following rule, which generalizes (F~2): $$N_{p+q}(x+y)=\sum_{\lambda\in
I}N_p(x+\lambda)\,N_q(y+\lambda^*)$$ for $x$, $y$ in ${\bf N}^{(I)}$, $p,q$ in
${\bf N}$
(compare with (4.5)).
\vskip0.7cm
6. {\it Diagonalization of the fusion rules}
\smallskip
 \ind To go further we need some information on the structure of the ring
${\cal F}$. We have
already observed  that ${\cal F}$ carries a symmetric, positive definite
bilinear form $<\ |\ >$
defined by $<\!x\,|\,y\!>= t(xy^*)$, for which  $I$ is an orthonormal basis.
The fact that
$^*$ is a ring homomorphism implies $<\!xy\,|\,z\!>=<\!x\,|\,y^*z\!>$  for all
$x,y,z$ in ${\cal F}$.
The existence of this form imposes strong restrictions on the ring ${\cal F}$.
 \medskip
{\pc PROPOSITION} 6.1$.-$ {\it The ${\bf Q}$\tx algebra ${\cal F}_{\bf
Q}:={\cal F}\otimes{\bf Q}$  is
isomorphic to a product $\prod K_i$ of finite extensions of ${\bf Q}$,
preserved by the involution,
which are of the following two types:}
\indp a) {\it a totally real extension of ${\bf Q}$ with the trivial
involution;}
\indp b) {\it a totally imaginary extension of ${\bf Q}$ which is a quadratic
extension of a
totally real extension  of ${\bf Q}$, the involution being the nontrivial
automorphism of
that quadratic extension.}
\ind (Recall that an
extension $K$ of ${\bf Q}$ is called totally real (resp. totally imaginary) if
$K\otimes_{\bf Q}{\bf
R}$ is isomorphic to ${\bf R}^{r_1}$ (resp. ${\bf C}^{r_2}$).) \ind We first
observe that the ring
${\cal F}$ is {\it reduced}: for $x\in{\cal F}$, the relation $x^2=0$ implies
$<\!xx^*\,|\,xx^*\!>=0$, hence $xx^*=0$, which in turn implies
$<\!x\,|\,x\!>=0$ and finally $x=0$.
Since a reduced finite-dimensional ${\bf Q}$\tx algebra is a product of fields,
we get the
decomposition ${\cal F}_{\bf Q}=\prod K_i$. This decomposition is canonical
(each factor corresponds
to an indecomposable idempotent of ${\cal F}_{\bf Q}$), so it is preserved by
the involution: each
factor $K_i$ is either preserved by $\sigma$, or mapped isomorphically onto
another factor $K_j$. In
the second case ${\cal F}_{\bf Q}$ contains a product of fields $K\times K$,
with the involution
interchanging the two factors. Then the set of elements $xx^*$ for $x\in
K\times K$ is the diagonal
$K\i K\times K$,  a ${\bf Q}$\tx  vector space on which $t$ can take arbitrary
values, contradicting
the positivity assumption. \ind Now let $K$ be one of the $K_i$'s, and  let
$\sigma$ denote the
induced involution.  Applying the same argument to the ${\bf R}$\tx algebra
$K\otimes_{\bf Q}{\bf R}$
we find that it is of the form ${\bf R}^{r_1}\times {\bf C}^{r_2}$, with
$\sigma$ preserving each
factor; since the induced involution on each factor is ${\bf R}$\tx linear,
there is no choice but
the identity on the real factors and the complex conjugation on the complex
ones. In particular, the
fixed subfield $K^\sigma$ of $\sigma$ is totally real and $K^\sigma\otimes_{\bf
Q}{\bf R}$ is
isomorphic to   ${\bf R}^{r_1+r_2}$. If we are not in case a), $K^\sigma$ is
strictly smaller than
$K$; counting degrees we get $$r_1+2r_2=[K:{\bf Q}]=2\,[K^\sigma:{\bf
Q}]=2(r_1+r_2)\ ,$$ hence
$r_1=0$, and we are in case b). \cqfd
\medskip
\ind Let $S$ be the set  of characters (i.e. algebra
homomorphisms) of ${\cal F}$ into ${\bf C}$; we can view $S$ as the spectrum of
the ${\bf C}$\tx
algebra ${\cal F}_{\bf C}:={\cal F}\otimes{\bf C}$. In the sequel we'll use
prop. 6.1 only through
the following  weaker corollary:
\smallskip
{\pc COROLLARY} 6.2$.-$ a) {\it The map ${\cal F}_{\bf C}\longrightarrow {\bf
C}^S$ given by
$x\mapsto (\chi(x))_{\chi\in S}$ is an isomorphism of ${\bf C}$\tx algebras.
\ind {\rm b)} One has $\chi(x^*)=\overline{\chi(x)}$ for} $\chi\in S$, $x\in
{\cal F}$.
\ind The assertion a) follows immediately from the proposition. We have seen in
the proof of the
proposition that ${\cal F}_{\bf R}$ is isomorphic as an algebra with involution
to  ${\bf
R}^p\times{\bf C}^q$, with the  involution acting trivially on the real factors
and by conjugation
on the complex factors; this is equivalent to b). \cqfd
\smallskip
\ind Clearly an explicit knowledge of the isomorphism ${\cal
F}_{\bf C}\longrightarrow {\bf C}^S$ (that is, of the characters $\chi:{\cal
F}\rightarrow {\bf C}$)
will allow us to perform any computation we need to in the ring ${\cal F}$. As
an example:
\smallskip
{\pc PROPOSITION} 6.3$.-$ {\it In the situation of prop.} 5.9, {\it one has }
$$N_g(\alpha_1+\ldots+\alpha_p)=\sum_{\chi\in
S}\chi(\alpha_1)\ldots\chi(\alpha_p)\,\chi(\omega)^{g-1}\qquad{with}\qquad
\chi(\omega)=\sum_{\lambda\in I}|\chi(\lambda)|^2\ .$$ \ind Let $x\in{\cal F}$;
the corresponding
element of ${\bf C}^S$ is $(\chi(x))_{\chi\in S}$. In the standard basis of
${\bf C}^S$, the matrix of
$m_x$ is the diagonal matrix with entries $(\chi(x))_{\chi\in S}$, so we have
$\Tr(x)=\sum_{\chi\in
S}\chi(x)$. Then the result follows from prop. 5.9. \cqfd
\medskip
\ind (6.4) One can obviously play for a while around these formulas; let me
give a sample, also to
make a link with the notation of the mathematical physicists.
 Let $\alpha\in I$; the matrix of the multiplication
$m_\alpha$ in the basis $I$ is
$N_\alpha=(N_{\alpha\gamma}^\beta)_{(\beta,\gamma)\in I\times I}$, with
$N_{\alpha\gamma}^\beta=N(\alpha+\beta^*+\gamma)$. On the other hand the matrix
of $m_\alpha$ in  the
standard basis of ${\bf C}^S$ is  the diagonal matrix $D_\alpha$ with entries
$\chi(\alpha)$, for $\chi\in S$. The base change matrix is
$\Sigma=(\chi(\lambda))_{(\chi,\lambda)\in S\times I}$, so that
$$N_\alpha=\Sigma^{-1}D_\alpha\Sigma$$ i.e. ``the matrix $\Sigma$ diagonalizes
the fusion rules''.
Observe that this remains true if we replace $\Sigma$ by
$\Sigma'=\Delta\Sigma$, where $\Delta$ is a
diagonal matrix. If we take  $\Delta=D_\omega^{-{1\over 2}}$ (noting that
$\chi(\omega)=\sum_\lambda
|\chi(\lambda)|^2$ is positive), an easy computation gives that the matrix
$\Sigma'$ is {\it
unitary}. This is only part of the story: for a RCFT the Verlinde conjecture
gives a geometric
interpretation of the matrix $\Sigma'$ in terms of the conformal blocks for
$g=1$, providing further
restrictions on the fusion ring ${\cal F}$.
\vfill\eject

\centerline{\bf Part III: The fusion ring ${\cal R}_{\ell}({\goth g})$}
\vglue 12pt plus 3pt 
7. {\it The rings ${\cal R}({\goth g})$ and
${\cal R}_{\ell}({\goth g})$} \ind (7.1) Recall that the representation ring
${\cal
R}({\goth g})$ is the Grothendieck ring  of finite-dimensional representations
of
${\goth g}$, with the multiplicative structure defined by the tensor product of
representations. It is a free ${\bf Z}$\tx module  with basis the isomorphism
classes of irreducible representations (i.e. the $[V_\lambda]$ for $\lambda\in
P_+$) with the rule
$$[V_\lambda]\cdot[V_\mu]=[V_\lambda\otimes V_\mu]\ .$$ \ind We are
interested in an analogue of ${\cal R}({\goth g})$ for the level
$\ell$ representations   of $\widehat{\goth g}$.
However it is not clear how to define  the multiplicative structure  in terms
of the
affine algebra $\widehat {\goth g}$: taking tensor products does not work,
since
the tensor product of two representations of level $\ell$ has level $2\ell$.
Instead
we will follow another route, which can be expressed purely in ordinary Lie
theory
terms.\medskip

\ind We have associated to the Lie algebra ${\goth g}$ and the integer
$\ell $ a fusion rule (example 5.2 {\it a}), defined by the formula $N(\sum
\lambda_i)=\dim V_{{\bf P}^1}(\vec P,\vec \lambda)$. We denote by ${\cal
R}_\ell({\goth g})$ the corresponding fusion ring, and call it the {\it fusion
ring of ${\goth g}$ at level} $\ell$. We can consider ${\cal R}_\ell({\goth
g})$ as the
free ${\bf Z}$\tx module with basis the isomorphism classes $[V_\lambda]$ for
$\lambda\in P_\ell$. The product in ${\cal R}_\ell({\goth g})$ is given by
$$[V_\lambda]\cdot[V_\mu]=\sum_{\nu\in P_\ell} N(\lambda+\mu+\nu^*)\,V_\nu\
$$ One can make this more explicit as follows. Consider the Lie subalgebra
${\goth
s}$ of ${\goth g}$ spanned by $H_\theta,X_\theta,X_{-\theta}$ (1.3); as in
(4.3)  we
denote  by $V^{(p)}$  the isotypic component of spin $p$ of a ${\goth s}$\tx
module
$V$.
 \smallskip
{\pc PROPOSITION} 7.2$.-$ {\it Let $\lambda,\mu\in P_\ell$. The product
$[V_\lambda]\cdot[V_\mu]$ in ${\cal R}_\ell({\goth g})$ is the class of the
${\goth g}$\tx module $V_\lambda\,{\scriptstyle\odot}\, V_\mu$ quotient of
$V_\lambda\otimes V_\mu$ by the ${\goth g}$\tx module spanned by the
isotypic components of spin $r$ of $V_\lambda^{(p)}\otimes V_\mu^{(q)}$ for all
triples $\{p,q,r\}$ such that $p+q+r>\ell$.}
 \ind By prop. 4.3, for each $\nu\in P_\ell$, $N(\lambda+\mu+\nu^*)$ is the
dimension of the space of ${\goth g}$\tx invariant linear forms  on
$V_\lambda\otimes V_\mu\otimes V_\nu^*$ which vanish on
$V_\lambda^{(p)}\otimes V_\mu^{(q)}\otimes (V_\nu^{(r)})^*$ for
$p+q+r>\ell$; this space is canonically isomorphic to  the space ${\cal
H}_{\lambda\mu}^\nu$ of ${\goth g}$\tx linear maps $u:V_\lambda\otimes
V_\mu\longrightarrow  V_\nu$ such that $\displaystyle
u(V_\lambda^{(p)}\otimes V_\mu^{(q)})\i \sum_{p+q+r\le \ell}V_\nu^{(r)}$,
that is such that $u$ annihilates the isotypic component of spin $r$ of
$V_\lambda^{(p)}\otimes V_\mu^{(q)}$ whenever  $p+q+r>\ell$.  Now for any
finite-dimensional ${\goth g}$\tx module $V$, the multiplicity of $V_\nu$ in
$V$
is $\dim\Hom_{\goth g}(V,V_\nu)$; therefore by definition  $\Hom_{\goth
g}(V_\lambda\,{\scriptstyle\odot}\, V_\mu,V_\nu)$ is isomorphic to the
subspace ${\cal H}_{\lambda\mu}^\nu$ of $\Hom_{\goth g}(V_\lambda\otimes
V_\mu,  V_\nu)$. This is equivalent by duality to the statement of the
proposition.
\cqfd  \smallskip {\it Examples} 7.3$.-$ {\it a}) Assume $\lambda+\mu\in
P_\ell$.
Then $V_\lambda^{(p)}$ and $V_\mu^{(q)}$ are nonzero only if $p\le{1\over 2}
\lambda(H_\theta)$ and $q\le {1\over 2}\mu(H_\theta)$, and
$V_\lambda^{(p)}\otimes
V_\mu^{(q)}$ has a component of spin $r$ if and only if $r\le p+q$.  Therefore
the condition $p+q+r\le \ell$ is always realized, so\break
$[V_\lambda]\cdot[V_\mu]$ {\it is the class of} $V_\lambda\otimes V_\mu$.
\ind {\it b}) Assume $\lambda(H_\theta)+\mu(H_\theta)=\ell+1$. Then  the
relation
$p+q+r> \ell$  holds if and only if $p={1\over 2}\lambda(H_\theta)\ ,\ q=
{1\over 2}\mu(H_\theta)\ ,\ r={1\over 2}(\ell+1)$; moreover every component
of spin ${1\over 2}(\ell+1)$ of $V_\lambda\otimes V_\mu$  occur in this way.
This
means that  $[V_\lambda]\cdot[V_\mu]$ {\it is obtained by removing from
$V_\lambda\otimes V_\mu$ all components $V_\nu$ with}
$\nu(H_\theta)=\ell+1$. \ind {\it c}) Assume
$\lambda(H_\theta)+\mu(H_\theta)=\ell+2$. Then  the same argument shows
that one has to remove the $V_\nu$'s with $\nu(H_\theta)=\ell+2\ {\rm or}\
\ell+1$, and the $V_\nu$'s with $\nu(H_\theta)=\ell$ which intersect
non-trivially
$V_\lambda^{(\lambda(H_\theta))}\otimes V_\mu^{(\mu(H_\theta))}$.
\vskip0.7cm
8. {\it The map ${\cal R}({\goth g})\longrightarrow {\cal R}_\ell({\goth g})$}
\ind Though ${\cal R}_\ell({\goth g})$ appears as a subgroup of ${\cal
R}({\goth
g})$, it is obviously not a subring. We will see, however, that there is a
natural way
to look at ${\cal R}_\ell({\goth g})$ as a {\it quotient ring} of ${\cal
R}({\goth
g})$.
\ind (8.1) We will need  a few classical facts about root systems, all of which
can be
found in [Bo].  For each root $\alpha$, the equation $\lambda(H_\alpha)=0$ (or
equivalently $(\lambda\,|\,\alpha)=0$) defines a hyperplane in  the real vector
space $P\otimes{\bf R}$, called the {\it wall} associated to $\alpha$. The {\it
chambers} of the root system are the connected components of the complement of
the walls. The chambers are fundamental domains for the action of the Weyl
group
$W$ on $P\otimes{\bf R}$.  \ind To the basis
$(\alpha_1,\ldots,\alpha_r)$ of the root system is associated a chamber $C$,
defined by the conditions $\lambda(H_{\alpha_i})\ge 0$. By definition the set
$P_+$
of dominant weights is $P\cap C$. Since $C$ is a fundamental domain, every
element
of $P$ can be written $w\lambda_+$ with $w\in W$, $\lambda_+\in P_+$; the
weight
$\lambda_+$ is uniquely determined, and so is $w$ if $\lambda$ does not belong
to
a wall. Let us   denote as usual by $\rho$ the half sum of the positive roots;
it is
characterized by the equality  $\rho(H_{\alpha_i})=1$ for each simple
root $\alpha_i$. Therefore  {\it the weights which belong to the interior of
$C$
are the weights $\lambda+\rho$ for $\lambda\in P_+$.} \ind For studying the
representation ring ${\cal R}_\ell ({\goth g})$ we need to consider  a closely
parallel situation where the role of $W$ is played by an infinite Coxeter
group,
the {\it  affine Weyl group} $W_\ell $. Let
\def\hv{h^{\scriptscriptstyle\vee}}
$\hv:=\rho(H_\theta)+1$\note{1}{This number is often called the {\it dual
Coxeter
number} of the root system.\medskip }. Then
$W_\ell $ is the group of motions of $P\otimes{\bf R}$ generated by $W$ and the
translation $x\mapsto x+(\ell +\hv)\theta$. Since each long root is conjugate
to $\theta$ under $W$, the group $W_\ell $ is the semi-direct product of $W$ by
the
lattice $(\ell +\hv)Q_{lg}$, where $Q_{lg}$ is the sublattice of $P$ spanned by
the
long roots. The {\it affine walls} of $P\otimes{\bf R}$ are the affine
hyperplanes
$(\lambda\,|\,\alpha)=(\ell +\hv)n$ for each root $\alpha$  and each $n\in{\bf
Z}$. The
connected components of the complement are called alcoves; each alcove is a
fundamental domain for the action of $W_\ell $ on $P\otimes{\bf R}$. The alcove
$A$
contained in $C$ and containing $0$ is defined by the inequalities
$\lambda(H_{\alpha_i})\ge 0$ for each basis root $\alpha_i$ and
$\lambda(H_\theta)\le \ell +\hv$.  We see as above that {\it the weights which
belong to the interior of $A$ are the $\lambda+\rho$ for $\lambda\in P_\ell$.}
 \ind Let ${\bf Z}[P]$ be the group ring of $P$; following [Bo] we denote by
$(e^\lambda)_{\lambda\in P}$ its canonical basis, so that the multiplication in
${\bf
Z}[P]$ obeys the usual rule $e^\lambda e^\mu=e^{\lambda+\mu}$. The action of
$W_\ell$ (hence of  $W$) on $P$ extends to an action on ${\bf Z}[P]$.
 Let $\varepsilon:W\rightarrow \{\pm 1\}$ be the signature homomorphism. We
denote by ${\bf Z}[P]_W$  the quotient of ${\bf Z}[P]$
by the sublattice spanned by the elements
$e^\lambda-\varepsilon(w)e^{w\lambda}$ $(\lambda\in P$, $w\in W)$ and by the
elements $e^\lambda$ for all weights $\lambda\in P$ belonging to a
wall\note{1}{Observe that such an element satisfies $\lambda(H_\alpha)=0$  for
some root $\alpha$, hence
$2e^\lambda=e^\lambda-\varepsilon(s_\alpha)e^{s_\alpha(\lambda)}$.};
we define   ${\bf Z}[P]_{W_\ell}$ in the same way. \medskip {\it Lemma} 8.2$.-$
{\it
The  linear maps}  $$\varphi:{\cal R}({\goth g})\longrightarrow {\bf
Z}[P]_W\qquad
,\qquad \varphi_\ell:{\cal R}_\ell({\goth g})\longrightarrow
 {\bf Z}[P]_{W_\ell}$${\it which associate to $[V_\lambda]$ the class of
$e^{\lambda+\rho}$, are bijective.}
\ind Let us define a linear map $\psi:{\bf Z}[P]\longrightarrow  {\cal
R}({\goth
g})$ in the following way: let $\lambda\in P$. By the above remarks, if
$\lambda$
does not lie on a wall,  there exist $w\in W$ and $\lambda_+\in P_+$, uniquely
determined,  such that
 $\lambda=w(\lambda_++\rho)$. We put $$\psi(e^\lambda)=\cases{
\varepsilon(w)\,[V_{\lambda_+}] & if $\lambda$  does not belong
to a wall,\cr \noalign{\smallskip} 0 & otherwise.\cr}$$ Then $\psi$ factors
through
$\overline{\psi}:{\bf Z}[P]_W\longrightarrow {\cal R}({\goth g})$, which
is easily seen to be the inverse of $\varphi$.  The same construction applies
 identically to define the inverse of $\varphi_\ell$. \cqfd
 \ind By the lemma there is a unique ${\bf Z}$\tx linear map $$\pi:{\cal
R}({\goth g})\longrightarrow {\cal R}_\ell({\goth g})$$  such that the diagram
$$\diagram{{\cal R}({\goth g})&\hfl{\pi}{}& {\cal R}_\ell({\goth g})\cr
\vfl{\varphi}{}&&\vfl{}{\varphi_\ell}\cr {\bf Z}[P]^{}_W&\hfl{p}{}&{\bf
Z}[P]^{}_{W_\ell} }$$where $p$ is the quotient map, is commutative. From the
lemma (and its proof) we get the following expression for $\pi$:
\medskip {\pc
PROPOSITION} 8.3$.-$  {\it  Let $\lambda\in P_+$; then \ind  --
$\pi([V_\lambda])=0$
if $\lambda+\rho$ belongs to an affine wall; \ind --
$\pi([V_\lambda])=\varepsilon(w)\,[V_\mu]$ otherwise, where $\mu\in P_\ell$,
$w\in W_\ell$ are such that $\lambda+\rho=w(\mu+\rho)$.  \ind In particular,
one
has $\pi([V_\lambda])=[V_\lambda]$ for $\lambda\in P_\ell$.} \cqfd
\vskip0.5cm
9. {\it  The spectrum of ${\cal R}_\ell({\goth g})$}
\smallskip
\ind (9.1) To understand the fusion ring ${\cal R}_\ell({\goth g})$ we need to
know
its  spectrum. Let us first consider the ring ${\cal R}({\goth g})$; it is
convenient to
introduce the simply-connected group $G$ whose Lie algebra is ${\goth g}$, and
the maximal torus $T\i G$ with Lie algebra ${\goth h}$. Any finite-dimensional
representation of ${\goth g}$  can be (and will be) considered as a $G$\tx
module. Any element $\lambda$ of $P$ defines a character $e^\lambda$ of $T$,
by the formula $e^\lambda(\exp H)=\exp \lambda(H)$; this defines an
isomorphism of $P$ onto the character group of $T$, which extends to an
isomorphism of the group algebra  ${\bf C}[P]$ onto the ring of algebraic
functions
on $T$.   We will identify  ${\bf Z}[P]$ to a subring of ${\bf C}[P]$, so that
the
notation $e^\lambda$ for the character associated to $\lambda$ is coherent with
the one we used before.

\ind (9.2) Each element $t$ of $T$ defines a character $\Tr_*(t)$ of ${\cal
R}({\goth
g})$, which associates to the class of a ${\goth g}$\tx module $V$ the number
$\Tr^{}_V(t)$. There is an explicit way of computing this character, the {\it
Weyl
formula}. Let me first introduce the antisymmetrization operator $J:{\bf
C}[P]\rightarrow {\bf C}[P]$, defined by the formula
$J(e^\mu)=\sum_{w\in W}\varepsilon(w)e^{w\mu}$; one has
$$J(e^\rho)=e^\rho\prod_{\alpha>0}(1-e^{-\alpha})$$([Bo], ch. VI, \S 3, prop.
2).
An  element  $t$ of $T$ is called {\it  regular} if $e^\alpha(t)\not=1$ for
each root
$\alpha$, or equivalently if   $wt\not=t$ for each  $w\in W$, $w\not=1$. Let
$t$ be a
regular element of $T$; one has $J(e^\rho)(t)\not=0$ and
$$\Tr^{}_{V_\lambda}(t)={J(e^{\lambda+\rho})(t)\over J(e^\rho)(t)}\ .$$
\medskip

\ind (9.3) We denote by  $T_\ell$  the subgroup of elements
$t\in T$ such that $e^{\alpha}(t)=1$
  for each element $\alpha$ of $(\ell +\hv)Q_{lg}$, and by $T^{\rm reg}_\ell $
the
subset of  regular  elements in $T_\ell $. The  finite group  $T_\ell$ will
play for
${\cal R}_\ell ({\goth g})$ the role of $T$ for ${\cal R}({\goth
g})$.
\smallskip {\it Lemma}$9.3.-$ a)  {\it For $t\in T^{\rm reg}_\ell $, the
character
$\Tr_*(t)$ factors through $\pi:{\cal R}({\goth g})\rightarrow {\cal
R}_\ell({\goth g})$.}  \ind b) {\it Let us
identify $P\otimes{\bf C}$ with  ${\goth h}$ using the normalized Killing form.
 Then the map $\displaystyle \lambda\mapsto \exp 2\pi
i{\lambda\over \ell +\hv}$ induces an isomorphism of $P/(\ell +\hv)Q_{lg}$
onto}
$T_\ell$.
\ind c) {\it The map $\displaystyle \lambda\mapsto \exp 2\pi
i{\lambda+\rho\over \ell +\hv}$ induces a bijection of $P_\ell $ onto}
$T_\ell^{\rm
reg}/W$.\smallskip

 \ind a) Let $t\in T^{\rm reg}_\ell $. The Weyl formula provides us with a
commutative diagram\vskip-15pt
$$\diagram{
{\cal
R}({\goth g})&&\cr \vfl{\varphi}{}&&\cr
{\bf Z}[P]_W&\hfl{}{j_t}&{\bf C}\cr
\fleche(10,18)\dir(3,-2)\long{15}
\put(20,16){\scriptstyle\Tr_*(t)}
}$$
\vskip-20pt where $j_t$ associates to the class of $e^\mu\in {\bf Z}[P]$ the
complex number $\displaystyle {J(e^\mu)(t)\over  J(e^\rho)(t)}$.
\smallskip
 \ind  The kernel of $\pi$ corresponds through
$\varphi$ to the kernel of $p$ (8.3), which is the subspace of ${\bf Z}[P]_W$
spanned
by the elements $e^{\mu+\alpha}-e^\mu$, for $\mu\in P$, $\alpha\in (\ell
+\hv)Q_{lg}$, and $e^\mu$ for $\mu$ in some affine wall. The elements of the
first type are killed by $j_t$ because $t$ is chosen so that $e^\alpha(t)=1$
for
$\alpha\in (\ell +\hv)Q_{lg}$; if $\mu$ belongs to an affine wall, $2e^\mu$ is
of the
first type (see the footnote to (8.1)), so one has  $2j_t(e^\mu)=0$ and
therefore
$j_t(e^\mu)=0$. This proves a).

\ind b)
Consider the exponential exact sequence

 $$0\rightarrow 2\pi
iQ^\vee\longrightarrow {\goth h}\ \hfl{\exp}{}\  T\rightarrow 0\ ;$$

here $Q^\vee$ is
the dual root lattice of the root system of ${\goth g}$, i.e. the lattice
spanned by
the $H_\alpha$'s. Let us denote by $P^\vee_{lg}$ the subgoup of
$Q^\vee\otimes{\bf
Q}$ consisting of elements $H$ such that $\alpha(H)\in{\bf Z}$ for all
$\alpha\in
Q_{lg}$; the map $\displaystyle H\mapsto \exp({2\pi i\over \ell +\hv}H)$
induces an
iso-\smallskip morphism of $P^\vee_{lg}/(\ell +\hv)Q^\vee$ onto $T_\ell $.
When
we identify $Q^\vee\otimes{\bf Q}$ with $P\otimes{\bf Q}$  using the normalized
Killing form,    $P^\vee_{lg}$ is identified with the dual lattice of $Q_{lg}$,
that is
the set of  elements $\lambda$ in   $P\otimes{\bf Q}$ such that
$(\lambda\,|\,\beta)\in {\bf Z}$ for each long root $\beta$. Because of the
normalization this is  equivalent to $\lambda(H_\beta)\in{\bf Z}$; since the
$H_\beta$'s   are the {\it short} roots of the dual system, and therefore span
the
coroot lattice $Q^\vee$, the dual lattice of $Q_{lg}$ is $P$.  In the same way
$Q^\vee$ is identified with the dual lattice $Q_{lg}$ of $P$. This proves b).
\ind c)
The isomorphism $P/(\ell +\hv)Q_{lg}\iso T_\ell $ is of
 course compatible with the action of $W$.
Now the orbits of $W$ in $P/(\ell +\hv)Q_{lg}$ are in one-to-one correspondence
with
the orbits of $W_\ell$ in $P$, and we have seen that those are parametrized by
the
elements of $P$ which lie in the affine alcove; moreover the orbits where $W$
acts
freely correspond to the weights which lie in the interior of the alcove, that
is
which are of the form $\lambda+\rho$ for $\lambda\in P_\ell.$ This gives c).
\cqfd
\medskip
 \ind (9.4) For $t\in T_\ell^{\rm reg}$, we will still denote by $\Tr_*(t)$ the
linear map ${\cal R}_\ell({\goth g})\longrightarrow {\bf C}$ obtained by
passing to
the quotient; because of prop.  8.3, it is again given by $[V]\mapsto
\Tr^{}_V(t)$. It
depends only on the class of $t$ in $T_\ell^{\rm reg} /W$. The two next results
are directly
borrowed from [F]:
\smallskip  {\pc PROPOSITION} 9.4$.-$ {\it The following
conditions are equivalent:} \ind (i) {\it The map $\pi:{\cal R}({\goth
g})\longrightarrow{\cal R}_\ell({\goth g}) $ is a ring homomorphism;}
\ind (ii) {\it One has $\pi([V_\lambda\otimes
V_\varpi])=[V_\lambda]\cdot[V_\varpi]$ for each $\lambda$ in $P_\ell$ and
each fundamental weight $\varpi\in P_\ell$};
 \ind (iii) {\it The linear forms} $\Tr_*(t)$ ($t\in T_\ell^{\rm reg} /W$) {\it
are
characters of the fusion ring ${\cal R}_\ell({\goth g})$.
\ind When these conditions hold, the spectrum  of ${\cal R}_\ell({\goth g})$
consists of the
characters $\Tr_*(t)$  where $t$ runs over $T_\ell^{\rm reg} /W$. }
\ind Since $\pi([V_\lambda])=[V_\lambda]$ for $\lambda\in P_\ell$ (prop. 8.3),
the
implication (i)$\Rightarrow $(ii) is clear.
\ind (ii)$\Rightarrow $(iii): Fix some $t\in T_\ell^{\rm reg} $, and put
$\chi=\Tr_*(t)$.
Let $\Omega$ denote the (finite) set of fundamental weights.
By the lemma below the ${\bf Z}$\tx algebra ${\cal R}_\ell({\goth g})$ is
generated by the family $([V_\varpi])_{\varpi\in\Omega}$. Therefore to prove
that $\chi$ is a character it is enough to check the equality
$\chi(x\cdot[V_\varpi]=\chi(x)\chi([V_\varpi])$ for  $x\in
{\cal R}_\ell({\goth g})$, $\varpi\in\Omega$; moreover because $\chi$ is ${\bf
Z}$\tx linear we may take $x$ of the form $[V_\lambda]$ for $\lambda\in
P_\ell$. But since
$\chi\rond\pi$ is a character of ${\cal R}({\goth g})$, this follows from (ii).
\ind  (iii)$\Rightarrow $(i):  Assume that (iii) holds. Different orbits of
$W$ in $T_\ell^{\rm reg} $ give different characters of ${\cal R}({\goth g})$,
hence of  ${\cal
R}_\ell({\goth g})$; because $ \Card(P_\ell)=\Card(T_\ell^{\rm reg} /W)$ by
lemma
9.3 {\it c}), the
spectrum  of ${\cal R}_\ell ({\goth g})$ consists of the characters  $\Tr_*(t)$
for
$t\in T_\ell^{\rm reg} /W$.  Since  $\Tr_*(t)\rond\pi$ is a ring homomorphism
for
each  $t$ in $T_\ell^{\rm reg} /W$ it follows that $\pi$ is a ring
homomorphism.~\cqfd
 {\it Lemma} 9.5$.-$ {\it The classes $[V_\varpi]$ for
$\varpi\in\Omega$  generate the ${\bf Z}$\tx algebra} ${\cal R}_\ell({\goth
g})$. \ind
Let us choose an element $H$ of $Q^\vee$ such that $\alpha_i(H)$ is a positive
integer for each simple root $\alpha_i\in B$. We will prove by induction on
$\lambda(H)$ that $[V_\lambda]$ is a polynomial in $([V_\varpi])_{\varpi\in
\Omega}$ for every $\lambda\in P_\ell$. This is clear if $\lambda=0$. If
$\lambda\not=0$, we can write $\lambda=\mu+\varpi$ with $\mu\in P_\ell$,
$\varpi\in\Omega$. Then  $V_\mu\otimes V_\varpi$ is the sum of $V_\lambda$ and
of irreducible ${\goth g}$\tx modules $V_{\nu}$  whose highest weights are of
the
form  $\nu=\lambda-\sum n_i\alpha_i$ with $n_i\in{\bf N}$, $\sum n_i>0$ (see
e.g.
[Bo], Ch. VIII, \S 7, \no 4, prop. 9).  By the induction hypothesis the
elements
$[V_\mu]$ and $[V_\nu]$ of ${\cal R}_\ell ({\goth g})$ are polynomial in the
$[V_\varpi]$'s; on the other hand the element $[V_\mu\otimes V_\varpi]$ is
equal
to $[V_\mu]\cdot[V_\varpi]$  (example 7.3 {\it a}). It follows that
$[V_\lambda]$
 is a   polynomial in the $[V_\varpi]$'s, hence the lemma. \cqfd

 {\pc PROPOSITION} 9.6 [F]$.-$ {\it The conditions of prop.  {\rm 9.4} hold
when
${\goth g}$ is of type $A_r$, $B_r$, $C_r$, $D_r$ or $G_2$.} \ind I will
content myself
with the cases $A_r$ and $C_r$, which are easy, and refer to the Appendix of
[F] for
the (rather technical) details in the remaining cases.  \ind By prop. 9.4 we
need to
prove the equality $\pi([V_\lambda\otimes
V_\varpi])=[V_\lambda]\cdot[V_\varpi]$ for each $\lambda$ in $P_\ell$ and each
fundamental weight $\varpi\in P_\ell$. For $A_r$ and $C_r$ a glance at the
tables in
[Bo] show that $\varpi(H_\theta)=1$ for each fundamental weight $\varpi$. If
$\lambda(H_\theta)<\ell$ we are done by example 7.3 {\it a}). If
$\lambda(H_\theta)=\ell$ we just have to apply the example 7.3 {\it b}) and
observe
that an irreducible ${\goth g}$\tx modules $V_\mu$ with
$\mu(H_\theta)=\ell+1$
is killed by $\pi$ (prop. 8.3). \cqfd

\smallskip
\ind We will now apply these results to reach our  goal, which is to compute
the dimension of the spaces  $V_C(\vec P,\vec \lambda)$ defined in part~I.
Since
we have now (at least in most cases) a precise description of the characters of
${\cal R}_\ell({\goth g})$, we can use the formula of prop.  6.3. The only
remaining
difficulty is to compute the expression $\sum_{\lambda\in P_\ell
}|\chi(\lambda)|^2$. This is provided by the following lemma:  \smallskip {\it
Lemma} 9.7$.-$ {\it Let $t\in T^{\rm reg}_\ell$. Then} $$\sum_{\lambda\in
P_\ell}|\Tr^{}_{V_\lambda}(t)|^2={|T_\ell |\over \Delta(t)}\ ,$$ where
$\displaystyle
\Delta(t):=|J(e^\rho)(t)|^2=\prod_{\alpha\in R({\goth g,h})}(e^\alpha(t)-1)$.
\ind For
$\lambda\in P_\ell$, let us denote by $t_\lambda$ the element  $\displaystyle
\exp
2\pi i {\lambda+\rho\over \ell +\hv}$ of $T$;  the  $t_\lambda$'s\smallskip
 for
$\lambda\in P_\ell$ form a   system  of representatives of $T^{\rm reg}_\ell/W$
(lemma 9.3 {\it c}). For $\lambda,\mu\in P_\ell$, one has
$$J(e^{\lambda+\rho})(t_\mu)=\sum_{w\in W}\varepsilon(w)\exp 2\pi
i{(w(\lambda+\rho)\,|\,\mu+\rho)\over \ell +\hv}=J(e^{\mu+\rho})(t_\lambda)\
$$
so by the Weyl formula (9.2)$$\sum_{\lambda\in
P_\ell}|\Tr^{}_{V_\lambda}(t_\mu)|^2={1\over
\Delta(t_\mu)}\sum_{\lambda\in P_\ell}|J(e^{\mu+\rho})(t_\lambda)|^2\
$$
  \ind Let $L^2(T_\ell)$ be the space of functions on the
(finite) group $T_\ell$, endowed with the usual scalar product
$$<\!f\,|\,g\!>={1\over |T_\ell|}\sum_{t\in T_\ell}\overline{f(t)}\,g(t)\ .$$
The function $h(t)=J(e^{\mu+\rho})(t)$ on $T$ is anti-invariant, i.e.
satisfies
$h(wt)=\varepsilon(w)h(t)$ for  $w\in W$, $t\in T$. It follows on one hand that
it
vanishes  at any non-regular point $t$ of $T$ (for any such point is fixed by a
reflection $s\in W$, so $h(t)=h(s(t))=-h(t)$), on the other hand that $|h|^2$
is $W$\tx
invariant. Therefore
$$\sum_{\lambda\in P_\ell}|J(e^{\mu+\rho})(t_\lambda)|^2={|T_\ell |\over
|W|}\ ||J(e^{\mu+\rho})||^2\ ,$$where the norm is taken in $L^2(T_\ell )$.
\ind  We claim that the restrictions to $T_\ell$ of the
characters $e^{w(\mu+\rho)}$, for $w\in W$, are all distincts. Suppose this is
not
the case; then there exists distinct elements  $w,w'\in W$ such that
$\bigl(w(\mu+\rho)-w'(\mu+\rho)\ |\ \lambda\bigr)\in (\ell +\hv){\bf Z}$ for
all $\lambda\in P$. But we have seen in the proof of lemma 9.3 {\it b}) that
the dual
lattice of $P$ is $Q_{lg}$, so the above condition means that
$\mu+\rho-w^{-1}w'(\mu+\rho)$ belongs to $(\ell +\hv)Q_{lg}$. This implies that
there exists a nontrivial element of $W_\ell $ fixing $\mu+\rho$, a
contradiction.
\ind Then  the orthogonality relations for the characters of the finite group
$T_\ell
$ give $||J(e^{\mu+\rho})||^2=|W|$, from which the lemma follows. \cqfd
\medskip
{\pc COROLLARY} 9.8 (Verlinde formula)$.-$ {\it
Assume  that the conditions of prop. {\rm 9.4} hold, e.g. that ${\goth g}$ is
of
type {\rm A, B, C, D} or {\rm G}. One has}
$$\eqalign{\dim V_C(\vec P,\vec \lambda)&=|T_\ell |^{g-1}\sum_{t\in T^{\rm
reg}_\ell } {\Tr^{}_{V_{\vec \lambda}}(t)\over \Delta(t)^{g-1}}\cr
&=|T_\ell |^{g-1}\sum_{\mu\in P_\ell } \Tr^{}_{V_{\vec \lambda}}(\exp
2\pi i {\mu+\rho\over \ell +\hv})
 \prod_{\alpha>0}\,\Bigl|2\sin{\pi{(\alpha\,|\,\mu+\rho)\over \ell
+h}}\Bigr|^{2-2g}\ .}$$
\ind The first expression is a simple reformulation of prop.  6.3 using the
explicit
description of the characters (prop. 9.4) and the above lemma. The second one
is obtained, after some easy manipulations, by using the description of
$T_\ell^{\rm reg}$ given in lemma 9.3 {\it c}). \cqfd \smallskip
{\it Remark} 9.9$.-$ One obtains easily an explicit expression for $|T_\ell |$,
using the isomorphism $P/(\ell +\hv)Q_{lg}\iso T_\ell $ (lemma  9.3 {\it b}).
One finds \  $|T_\ell |=(\ell +\hv)^rfq$,
\ where $r$ is the rank of ${\goth g}$, $f$ its connection index $(=|P/Q|)$,
and
$q$ the index of $Q_{lg}$ in $Q$. A glance at the tables gives $q=2$ for $B_r$,
$2^{r-1}$ for $C_r$, $4$ for $F_4$, $6$ for $G_2$, and of course $1$
otherwise.
 \vfill\eject \centerline{ REFERENCES} \vskip0.7cm \baselineskip14pt
\def\num#1{\item{\hbox to\parindent{\enskip [#1]\hfill}}}
\parindent=1.5cm

\num{B-L} A. {\pc BEAUVILLE}, Y. {\pc LASZLO}: {\sl
Conformal blocks and generalized theta functions.} Comm.  Math.
Phys., to appear.
\smallskip

\num{Bo} N. {\pc BOURBAKI}: {\sl Groupes et Alg\`ebres de Lie}, ch. 4 to 8.
Masson,
Paris (1990).\smallskip

 \num {F } G. {\pc FALTINGS}: {\sl A proof for the Verlinde formula.} J.
Algebraic Geometry, to appear. \smallskip
\num K V. {\pc KAC}: {\it Infinite dimensional Lie algebras} (3rd edition).
Cambridge University Press (1990). \smallskip

\num{S} A. {\pc SZENES}: {\sl The combinatorics of the Verlinde formula.}
Preprint (1994).\smallskip
  \num{T-U-Y} A. {\pc TSUCHIYA}, K. {\pc
UENO}, Y. {\pc YAMADA}: {\sl  Conformal field theory on universal family of
stable
curves with gauge symmetries.} Adv. Studies in Pure Math. {\bf 19}, 459-566
(1989).
\smallskip
 \num V E. {\pc VERLINDE}: {\sl Fusion rules and modular transformations
in 2d conformal field theory.} Nuclear Physics {\bf B300}  360-376
(1988). \smallskip
\vskip1cm
\hfill\hbox to 5cm{\hfill A. Beauville\hfill}

\hfill\hbox to 5cm{\hfill URA 752 du CNRS\hfill}

\hfill\hbox to 5cm{\hfill Math\'ematiques -- B\^at. 425\hfill}

\hfill\hbox to 5cm{\hfill Universit\'e Paris-Sud\hfill}

\hfill\hbox to 5cm{\hfill 91 405 {\pc ORSAY} Cedex, France\hfill}
\bye